\documentclass[11pt,letterpaper,final]{article}
\usepackage[utf8]{inputenc}
\usepackage[margin=0.96in]{geometry}
\usepackage{enumerate}
\usepackage{amsmath}
\usepackage{amsfonts}
\usepackage{amssymb}
\usepackage{graphicx}
\usepackage[font={it}]{caption}
\usepackage{color, float}
\usepackage[dvipsnames]{xcolor}
\usepackage{hyperref}
\hypersetup{
   colorlinks=true,
   linktoc=all,
   linkcolor=blue,
}
\usepackage{listings}
\usepackage{xcolor}

\definecolor{codegreen}{rgb}{0,0.6,0}
\definecolor{codegray}{rgb}{0.5,0.5,0.5}
\definecolor{codepurple}{rgb}{0.58,0,0.82}
\definecolor{backcolour}{rgb}{0.95,0.95,0.92}

\lstdefinestyle{mystyle}{
    backgroundcolor=\color{backcolour},   
    emph={U_skinny, U_Skinny_to_U},
    emphstyle={\color{blue}},
    commentstyle=\color{purple},
    keywordstyle={\color{teal}},
    numberstyle=\tiny\color{violet},
    stringstyle=\color{codepurple},
    basicstyle=\ttfamily\footnotesize,
    breakatwhitespace=false,         
    breaklines=true,                 
    captionpos=b,                    
    keepspaces=true,                 
    numbers=left,                    
    numbersep=5pt,                  
    showspaces=false,                
    showstringspaces=false,
    showtabs=false,                  
    tabsize=2,
}

\title{\bf Techniques, Tricks and Algorithms for Efficient GPU-Based Processing of \\
Higher Order Hyperbolic PDEs}
\author{Sethupathy Subramanian$^1$ {\small (ssubrama@nd.edu)},  
        Dinshaw S. Balsara$^{1,2,*}$ {\small (dbalsara@nd.edu)}, \\
        Deepak Bhoriya$^{1,3}$ {\small (dbhoriy2@nd.edu)} and 
        Harish Kumar$^3$ {\small (hkumar@iitd.ac.in)}}
\date{%
    {\small \it
    $^1$Physics Department, University of Notre Dame, IN, USA, 
    $^2$ACMS Department, University of Notre Dame, IN, USA, 
    $^3$Mathematics Department, Indian Institute of Technology Delhi, India, 
    $^*$Corresponding author} 
    }

\renewcommand{\thesection}{\Roman{section}}

\begin{document}
\maketitle
\tableofcontents

\section*{Abstract}

    GPU computing is expected to play an integral part in all modern Exascale
supercomputers. It is also expected that higher order Godunov schemes will make up about a
significant fraction of the application mix on such supercomputers. It is, therefore, very
important to prepare the community of users of higher order schemes for hyperbolic PDEs for
this emerging opportunity.

    Not every algorithm that is used in the space-time update of the solution of hyperbolic
PDEs will take well to GPUs. However, we identify a small core of algorithms that take
exceptionally well to GPU computing. Based on an analysis of available options, we have been
able to identify Weighted Essentially Non-Oscillatory (WENO) algorithms for spatial
reconstruction along with Arbitrary DERivative (ADER) algorithms for time-extension followed
by a corrector step as the winning three-part algorithmic combination. Even when a winning
subset of algorithms has been identified, it is not clear that they will port seamlessly to GPUs.
The low data throughput between CPU and GPU, as well as the very small cache sizes on
modern GPUs, implies that we have to think through all aspects of the task of porting an
application to GPUs. For that reason, this paper identifies the techniques and tricks needed for
making a successful port of this very useful class of higher order algorithms to GPUs.

    Application codes face a further challenge the GPU results need to be practically
indistinguishable from the CPU results in order for the legacy knowledge-bases embedded in
these applications codes to be preserved during the port of GPUs. This requirement often makes
a complete code rewrite impossible. For that reason, it is safest to use an approach based on
OpenACC directives, so that most of the code remains intact (as long as it was originally well-
written). This paper is intended to be a one-stop-shop for anyone seeking to make an OpenACC-
based port of a higher order Godunov scheme to GPUs.

    We focus on three broad and high-impact areas where higher order Godunov schemes are
used. The first area is computational fluid dynamics (CFD). The second is computational
magnetohydrodynamics (MHD) which has an involution constraint that has to be mimetically
preserved. The third is computational electrodynamics (CED) which has involution constraints
and also extremely stiff source terms. Together, these three diverse uses of higher order Godunov
methodology, cover many of the most important applications areas. In all three cases, we show
that the optimal use of algorithms, techniques and tricks, along with the use of OpenACC, yields
superlative speedups on GPUs! As a bonus, we find a most remarkable and desirable result:
some higher order schemes, with their larger operations count per zone, show better speedup
than lower order schemes on GPUs. In other words, the GPU is an optimal stratagem for
overcoming the higher computational complexities of higher order schemes! Several avenues for
future improvement have also been identified. A scalability study is presented for a real-world
application using GPUs and comparable numbers of high-end multicore CPUs. It is found that
GPUs offer a substantial performance benefit over comparable number of CPUs, especially when
all the methods designed in this paper are used.

{\bf Mathematics Subject Classification:} 65M08, 65M22, 76M12, 35Q35, 35Q61, 35Q68, 35Q85, 76W05, 68Q85

{\bf Key Words:} PDEs, numerical schemes, mimetic, high performance computing

\section{Introduction}
    The advent of GPUs with floating point support and error correction in hardware, as well
as the advent of open standards for programming GPUs with OpenACC directives, could
produce a sea change in Exascale computing. Almost every major high performance computing
(HPC) unit is on track to install an Exascale supercomputer in the next few years. Such
supercomputers will consist of many nodes. Each node will have a multicore CPU feeding data
to several on-node GPUs. It is expected that applications codes that use higher order Godunov-
based algorithms will make up a significant fraction of the use cases for such Exascale
supercomputers. As a result, it is very important, and very timely, to make a deep analysis of
different algorithms that are available under the rubric of higher order Godunov (HOG) schemes.
Some of these algorithms take well to GPU computing, but others can be seen from the very
onset to be weaker options. Even with a favorable algorithmic mix, it is realized that GPUs do
have their own inherent weaknesses and limitations. Even for the winning set of algorithms, it is
very important to identify those techniques and tricks that allow the full potential of the GPU to
be realized while skirting around the weaknesses of the GPU! The goals of this paper are,
therefore, fourfold:
\begin{itemize}
    \item Identify the winning combination of algorithms that result in a HOG scheme that will port
successfully to GPUs.
    \item Identify the available set of techniques and tricks that allow these HOG algorithms to be
(painlessly) ported over for GPU computing. This is done within the context of OpenACC so that
the interested reader should be able to pick up the minimal set of directives needed for porting a
HOG code to a GPU just by reading this one article.
    \item Show that this salubrious combination of optimal algorithm sets and the right set of techniques
and tricks yields very attractive performance gains on modern GPUs.
    \item Show that these advantages extend to different hyperbolic systems with different collocations
and constraints and treatment of stiff source terms. This gives diverse communities of scientists
and engineers confidence that the innovations reported here are broad-based and applicable to
their own community.
\end{itemize}

To give the reader a foretaste of what is to come, we will show that as one goes to HOG schemes
with increasing order of accuracy, the speedup in going from CPU to GPU actually improves
(with some caveats)! In other words (and under favorable circumstances), the GPU is an ideal
vehicle for offsetting the higher cost of an algorithm with greater computational complexity!
    In this paper we will focus on three different codes from three different application areas,
all of which have been ported to GPUs. These codes represent the different styles in which HOG
schemes are used in many prominent applications areas. The three focus areas, and our
justification for choosing them, are:-

\begin{enumerate}
\item We will focus on a code drawn from computational fluid dynamics (CFD). Such codes are
workhorses in numerous engineering applications and also in aerospace and geophysical
modeling arenas. CFD codes solve the compressible Euler equations and rely on a conventional
HOG structure with zone-centered fluid variables that are updated with face-centered fluxes.
CFD codes exemplify a traditional application of classical HOG schemes.

\item The second code that we focus on is a code drawn from magnetohydrodynamics (MHD). Such
codes contribute to applications in astrophysics, space physics, plasma physics, fusion and space
re-entry vehicles. They rely on a solution of the compressible MHD equations, which have
extremely non-linear and non-convex terms that can give rise to compound waves. The MHD
equations have an inbuilt divergence constraint; i.e., the divergence of the magnetic field has to
remain zero. This requires a face-centered collocation of the magnetic fields which are updated
with edge-centered electric fields. As a result, an MHD code retains fluid variables that are
collocated in the same way as a CFD code along with a divergence constraint-free evolution of
the face-centered magnetic fields. Therefore, MHD codes exemplify a class of problems with
extremely strong non-linearities and mimetic constraint-preservation.

\item The third code that we focus on is a code drawn from computational electrodynamics (CED).
Such codes are very important in electrical and electronic engineering. Example scientific and
engineering applications include: understanding the effects (and mitigation thereof) of solar
storms upon power grids;  further reductions of the radar cross-section of military platforms;
improvements in wireless communications and sensing by incorporating millimeter wave and
terahertz technologies; sophisticated first-principles analysis and design of nanoscale
plasmonics-based technologies for ultrafast optical computers and ultra-efficient photovoltaics;
and even the development of novel optical microscopy techniques to reliably detect deadly
human cancers at the earliest possible stage. CED codes evolve a combination of Faradays law
and the extended Amperes law using modern upwind methods and constraint-preserving
reconstruction. The magnetic induction and the electric displacement are both facially-collocated
and are updated by using edge-collocated electric and magnetic fields. Besides the preservation
of divergence constraints for both the facially collocated vector fields, such CED codes have to
handle extremely stiff source terms because the conductivity of some materials can take on
exceptionally large values. Therefore, CED codes exemplify a class of problems with mimetic
constraint-preservation and extremely stiff source terms.
All three codes in our application mix use ideas drawn from finite volume-type weighted
essentially non-oscillatory (WENO) reconstruction that are adapted to the different collocations
that are used. For each case, we will demonstrate the functioning of second, third and fourth
order schemes so that we can demonstrate sustained GPU speedup at several orders of accuracy.
In the next three paragraphs we discuss each of these codes with a view to identifying the
algorithmic sub-set that is best suited for GPU computing. When reading the next three
paragraphs, the reader should realize that selecting GPU-friendly algorithms is like threading the
needle! Out of the many algorithmic possibilities for designing a CFD, MHD or CED code, only
a few satisfy all the requirements that we need for obtaining superlative performance on a GPU.
\end{enumerate}

\subsection{Identifying Optimal Algorithmic Combinations for GPU-based Processing of HOG Schemes}
    CFD was indeed the first application area for which HOG schemes were developed
(Godunov [27], van Leer [42], Woodward and Colella [43], Harten et al. [29], Shu and Osher
[34], [35]). Early methods for CFD resorted to second order reconstruction of zone-centered flow
variables. In time, several authors (Jiang and Shu [30], Balsara and Shu [9], Balsara et al. [10],
Balsara, Garain and Shu [12]) provided very efficient formulations for implementing WENO
reconstruction at all orders on structured meshes. When these reconstruction strategies are
implemented directly for the conserved variables, the code can be written to absolutely minimize
cache usage, which makes modern WENO reconstruction strategies one of the favored
algorithmic elements for GPU computation. Progress was also made in the design of Riemann
solvers (Rusanov [32], van Leer [42], Colella [17], Roe [31], Harten, Lax and van Leer [28],
Einfeldt et al. [25], Toro, Spruce and Speares [40], Dumbser and Balsara [24]), to the point
where modern Riemann solvers have a very low operation count and use cache memory very
sparingly. Such Riemann solvers are to be favored in GPU environments because each core of a
GPU has a very small cache. In particular, the Rusanov Riemann solver (Rusanov [32]), the HLL
Riemann solver (Harten, Lax and van Leer [28]), the HLLC Riemann solver (Toro, Spruce and
Speares [40]) and the HLLI Riemann solver (Dumbser and Balsara [24]) have very favorable
cache usage.

    It was realized early on that just the zone-centered reconstruction, coupled with
application of Riemann solvers at the zone faces, would yield a very effective single stage of a
multi-stage, high order in time, Runge-Kutta timestepping scheme (Shu and Osher [34], Shu
[36]). But please realize that in Exascale computation, patches of the mesh will be assigned to
different nodes. Each patch will have to transfer boundary data to other patches that reside on
other nodes. Consequently, Runge-Kutta schemes would require back and forth data transfers at
each and every stage from host to device and vice versa, which is why they would be somewhat
disfavored on GPUs. (Due to their immense popularity, we do nevertheless include Runge-Kutta
schemes in the speed comparisons that we present in the results section.) Arbitrary DERivatives
in space and time (ADER) schemes for higher order timestepping were designed based on
Cauchy-Kovalevskaya methods (Titarev and Toro [38], [39], Toro and Titarev [41]) but those
are again very memory intensive at higher orders, resulting in them being disfavored for GPU
usage. In addition, for general hyperbolic systems, the Cauchy-Kovalevskaya procedure is
unusually difficult to formulate for larger hyperbolic systems, with the result that it is not a
general-purpose method that will extend to all manner of hyperbolic PDEs. More recent
incarnations of ADER schemes (Dumbser et al. [22], Balsara et al. [10]) are more suited to low
cache usage and also generalize well to several different PDEs. For that reason, we will favor
ADER timestepping schemes. To summarize the results of this paragraph, we arrive at a strategy
for CFD that first relies on finite volume WENO reconstruction; followed by an ADER predictor
step; and completed by a corrector step that invokes certain suitable Riemann solvers. Such a
strategy, which is based on three conceptual stages WENO reconstruction, ADER predictor and
Riemann solver-based corrector  was found to be an optimal approach for GPU computing.
(ADER-WENO methods are extensively discussed in a review by Balsara [8].) In this approach,
there is only one data transfer from host to node at the start of a timestep followed by a final data
transfer at the end of a timestep from node to host. Even this data transfer will be minimized, as
we will see in Section II.2. Also please note that due to the popularity and proliferation of
Runge-Kutta based time stepping, this option will also be explored in this paper.

    In addition to the fluid variables, the MHD equations include Faradays law for the
evolution of the divergence-free magnetic field. Mimetic HOG schemes that preserve the
divergence constraint have been designed (Ryu et al. [33], Dai and Woodward [21], Balsara and
Spicer [1]). Such schemes collocate the normal components of the magnetic field vector at the
faces of the mesh, and they are updated using the components of the electric field vector that is
collocated at the edges of the mesh. Progress was made on the second order reconstruction of the
magnetic field in Balsara [2], [3] and the same was extended to higher orders using WENO-like
methods in Balsara [4]. These WENO-like methods retain the advantage of minimizing cache
usage, making them optimal for GPUs. ADER methods were first applied to MHD in Balsara et
al. [10], [11] and these ADER methods also retain the same advantages as before. Evaluation of
the edge-centered electric fields requires multidimensional upwinding at the edges and that was
achieved by HLL type multidimensional Riemann solvers (Balsara [5]) along with their HLLC
variants (Balsara [6]) and their HLLI-type extensions (Balsara [7], Balsara and Nkonga [14]).
These multidimensional HLL-style Riemann solvers retain the same advantages of cache
friendliness as their one-dimensional versions, making them optimal for use on GPUs. Therefore,
we see that all the algorithmic pieces that were identified for CFD find their extensions for
MHD, making it possible to design optimal algorithmic combinations for use in GPU-enabled
MHD codes.

    CED shares some of the mimetic features with MHD in that one is now evolving two
curl-type equations; one for Faradays law and the other for the extended Ampere law. Face-
centered collocations are, therefore used for the normal component of the magnetic induction and
electric displacement. These two vector fields are then updated using edge-centered components
of the electric field and the magnetic field. WENO-like reconstruction strategies for the face-
centered primal variables were documented in Balsara et al. [15], [16]. Multidimensional
Riemann solvers of the HLL type for CED were also documented in Balsara et al. [15]. The
CED equations incorporate stiff source terms by way of conductivities that can sometimes be
millions of times larger than the contribution from the wave-propagation terms. As a result, an
optimal ADER formulation was documented in Balsara et al. [15], [16]. The ADER formulation
is special because it allows the stiff source terms to be treated implicitly while treating the rest of
the hyperbolic system explicitly. The implicit treatment requires the inversion of modest-sized
matrices locally within each zone. Furthermore, the unique ADER formulation mentioned above
decomposes the implicit treatment into a set of smaller block matrices whose inversion can be
carried out independently. As a result, the matrix inversion is very easy, very tractable, though
not very cache friendly. Consequently, we see that all the algorithmic pieces that were identified
for CFD find their extensions for CED, making it possible to design optimal algorithmic
combinations for use in GPU-enabled CED codes.

    The above three paragraphs have shown how one (figuratively) can thread the
algorithmic needle, arriving at end-to-end algorithmic ingredients for HOG schemes that are
ideally suited for GPU processing. In Section II we describe how these are assembled. Section III
shows results from GPU speedups for CFD, MHD and CED codes on A100 GPUs. Section IV
offers some conclusions. Since the examples in Section II used Fortran-based pseudo code to
exemplify the structure of GPU-friendly HOG schemes, in the Appendix we provide their C++
equivalents. This gives practitioners a one-stop-shop from which they can pick up the essentials
of GPU programming for HOG schemes. Section IV shows that the accuracy of the code is
virtually unaffected by whether it is run on a GPU or a CPU; thereby giving us confidence that
GPU computation should give us results that are entirely comparable in quality to the ones that
we have been used to getting on CPUs. Section V takes us through the process of setting up a
real-world simulation on a GPU-rich supercomputer and carrying out a scalability study on such
a machine. The exercise is intended to give us insight into the limitations as well as the
opportunities in GPU computing. Section VI presents some conclusions.

\section{GPU Implementation}
\subsection{GPU Architecture and Programming Environment and the Resulting Constraints on
Efficient Implementation}
    The trend in Exascale computing is to compute with the lowest power consumption,
because the cost of keeping an Exascale supercomputer running can indeed be quite prohibitive.
A majority of the end-to-end cost of owning and maintaining a supercomputer over its life cycle
goes towards the cost of electricity. GPUs are the preferred pathway for Exascale computing
because they couple performance with a level of (low) power consumption that cannot be
matched by modern CPUs. As a result, a modern Exascale supercomputer will consist of many
nodes that are connected together with a fast interconnect. Each such node will have a multicore
CPU that is coupled with multiple high-end GPUs. The V100 and A100 GPUs from Nvidia are
good exemplars of such GPUs; though at the time of this writing AMD has also announced a
GPU product. The compute power offered by one of those GPUs easily exceeds the compute
power offered by a multicore CPU. The possibility of having 4 or 6 such GPUs per node means
that the majority of on-node computation will be done by GPUs. It is for this reason that it is
beneficial to very briefly discuss the architecture of a modern GPU, especially as it pertains to
higher order Godunov schemes. Applications codes that use such schemes are anticipated to
constitute a significant fraction of the workload on Exascale supercomputers, making the
problem of optimal GPU performance very useful and timely.

    The node is often referred to as the host because it orchestrates the majority of the data
structure management on an Exascale supercomputer. The GPUs are referred to as the devices.
There is a bus that connects the host to the node, but that bus is by necessity quite slow. As a
result, the bus can act as a severe bottleneck and maintaining optimal performance requires
managing data. Modern GPUs that are used for HPC can be thought of as extensive arrays of
rather slow and weak processors (known as cores) with very limited cache memory. As a result,
cache has also to be used very efficiently. The fact that each GPU core functions as a slow and
weak processor is not to be scoffed at, because a modern GPU can have as many as 5000 such
cores. Taken together, all the cores of a GPU can offer performance that is superior to a modern
multicore CPU. What has made modern GPUs very desirable for HPC is the fact that the GPU
cores also offer double precision floating point performance as well as error correction in
hardware with the result that it is possible to expect that a result obtained on a GPU can be as
accurate as a result obtained on a CPU.

    It should also be mentioned that the amount of memory available on present-day GPUs is
often capped off at ~ 80GB, with the result that GPUs have much smaller integrated memory
than the total RAM memory on a CPU node. As a result, an Exascale computation will have to
be chunked into small chunks which are sent to the GPU for processing. For a structured mesh
application, which we address here, this means that the node may hold multiple smaller chunks
(or patches) of data, each of which is sent to the GPU for a time-update. Because each patch will
have to provide ghost zone information to neighboring patches, it is best to keep those patches
resident on the host and have the host manage the MPI-based communication. Optimal one-sided
communication strategies have been discussed in Garain, Balsara and Reid [26] and will not be
discussed here. The focus of our study in this paper is to optimize on-node GPU usage for higher
order Godunov schemes.

    Modern higher order Godunov schemes have a very intricate algorithmic structure and
codes that encapsulate such algorithms can be many thousands of lines long. Some such codes
can represent decades worth of effort and they need to reproduce certain well-known and well-
calibrated benchmark computations. A complete code rewrite, for example in CUDA, is often
unfeasible mainly because it would destroy the codes ability to reproduce the well-known and
well-calibrated benchmarks. As a result, it is beneficial to find the quickest pathway to port those
codes to GPUs with minimal rewriting of code. The advent of reasonably mature compilers that
support OpenACC gives us good grounds for thinking that this can be accomplished. OpenACC
is not a language, but rather, a set of language extensions. Similar to OpenMP, it has a set of
directives that can direct the computer to express the parallelism that is inherent in an algorithm.
Textbooks that describe OpenACC are available (Chandrasekaran and Juckeland [18]). For the
sake of completeness, it should also be mentioned that OpenMP is also expected to have GPU
support and textbooks for OpenMP include (Chapman, Jost and van der Pas [19]). Indeed, a code
that is well-structured for OpenMP can be very quickly adapted to OpenACC. It is our
impression that OpenACC currently gives more fine-grained control over GPU computations and
that the supporting compiler technology is also more mature. For this reason, we focus attention
on GPU computation of higher order schemes using OpenACC directives. This is the topic for
the rest of this section. The reader who is familiar with OpenMP will see many parallels between
OpenMP and OpenACC; and an OpenMP code can be an excellent starting point for a port to
OpenACC.

\subsection{Overall Structure of a Higher Order Godunov Scheme Implementation and its
Adaptation to GPUs}

    All HOG schemes of the ADER-WENO type have a conceptually identical structure;
even if the different collocation of variables require them to be different in detail. All such
schemes begin with a spatial reconstruction of the primal variables with higher order accuracy.
This first step examines the neighboring primal variables to build a higher order representation of
the primal variables within each zone. This higher order representation is provided by the
WENO or WENO-like reconstruction step. Once we know the spatial variation of all the
variables in a hyperbolic PDE, the time-dependent structure of the PDE can itself be exploited to
give us an in-the-small evolution of these variables in the temporal direction. This temporal
evolution has to be sufficiently accurate so as to match the accuracy of the spatial reconstruction.
This second step is provided by the ADER predictor step. The third step consists of a Riemann
solver-based corrector step. In other words, the space-time evolved (predicted) solution of the
hyperbolic PDE within each zone is made to interact with the analogous solutions from
neighboring zones. This interaction is mediated by Riemann solvers which mediate the
resolution of any jumps in the solution in a fashion that respects the physics of the PDE. This
corrector step yields the numerical fluxes that are used for a space and time accurate update of
the primal variables of the PDE. Because all ADER-WENO type algorithms for CFD, MHD and
CED share this same three-part structure, we will illustrate our ideas by instantiating them within
the context of a second order accurate CFD scheme.

\begin{figure}
\begin{lstlisting}[language=Fortran]
PROGRAM PSEUDO_HYDRO

INTEGER, PARAMETER :: nx = 48, ny = 48, nz = 48
REAL U(-1:nx+2, -1:ny+2, -1:nz+2, 5, 5)
REAL U_skinny(-1:nx+2, -1:ny+2, -1:nz+2, 5), dU_dt(nx, ny, nz, 5)
REAL, DIMENSION(0:nx, 0:ny, 0:nz,5) :: Flux_X, Flux_Y, Flux_Z
REAL dt, dt_next, cfl, dx, dy, dz

number_of_timesteps = 500

!$acc data create( U_skinny, U, dU_dt, Flux_X, Flux_Y, Flux_Z,
!$acc&        dt, dt_next, cfl, dx, dy, dz )

CALL Initialize(nx, ny, nz, U_skinny, dt, cfl, dx, dy, dz)  ! CPU subroutine

!$acc update device( cfl, dx, dy, dz )

DO i_timestep = 1, number_of_timesteps

!$acc update device( U_skinny, dt )
  CALL U_Skinny_to_U( nx, ny, nz, U_skinny, U )
  CALL Boundary_Conditions( nx, ny, nz, U )
  CALL Limiter( nx, ny, nz, U )
  CALL Predictor( nx, ny, nz, U )
  CALL Make_Flux_X( nx, n y, nz, U, Flux_X )
  CALL Make_Flux_Y( nx, ny, nz, U, Flux_Y )
  CALL Make_Flux_Z( nx, ny, nz, U, Flux_Z )
  CALL Make_dU_dt ( nx, ny, nz, Flux_X, Flux_Y, Flux_Z, dU_dt,
                    dt, dx, dy, dz )
  CALL Update_U_Timestep( nx, ny, nz, U, U_skinny, dU_dt, dt_next,
                          cfl, dx, dy, dz )

!$acc update host( U_skinny, dt_next )
  dt = dt_next
END DO

!$acc end data
END PROGRAM PSEUDO_HYDRO
\end{lstlisting}
    \caption{Showing the structure of a hydro code (in pseudocode format) at second order with
OpenACC extensions to allow it to perform well on a GPU. Notice that most of the code is identical
    to a serial code where we have a variable {\tt U} which holds the five fluid components and their
    space-time variation. The data at the start of a timestep is located in {\tt U( :, :, :, 1:5, 1)}. 
    This data is used to obtain the spatial variation which is stored in {\tt U( :, :, :, 1:5, 2:4)} 
    by using the Limiter routine. The Predictor routine builds the temporal dependence in 
    {\tt U( :, :, :, 1:5, 5)}. The {\tt Make\_Flux} routines make the appropriate fluxes and these 
    are used in {\tt Make\_dU\_dt} to make the time rate of change which is used in Update to make the time 
    update. These are all standard components of a classical hydro code.
    The changes that we introduced are shown using color. The blue variable {\tt U\_skinny} is the only
minimal amount of data that is moved from host to device at the beginning of the timestep and then
back from device to host at the end of a timestep. The OpenACC pragmas for carrying out this
data motion are shown in red. An OpenACC pragma is also needed to tell the device to create
requisite data on the device. The only new subroutine that we had to introduce is shown in
    blue, called {\tt U\_Skinny\_to\_U}. It operates on the device and takes the {\tt U\_Skinny} variable
    and copies it into {\tt U( :, :, :, :, 1)} at the beginning of each timestep}
\end{figure}

    The general structure of the HOG scheme is shown in Fig. 1, along with the necessary
OpenACC directives. (While we use Fortran in the main body of the text, the figures in the
Appendix show the exact C++ equivalents, so that both major programming paradigms are
covered here.) The conserved variable {\tt U} is dimensioned as {\tt U(-1:nx+2, -1:ny+2, -1:nz+2, 5, 5)}, 
here the first three dimensions correspond to the zones along the x, y, z directions
respectively. The active zones are numbered from {\tt 1:nx, 1:ny, 1:nz} and there are two ghost
cells on the either side, which are necessary for a second order accurate HOG scheme. The next
dimension 5 corresponds to the 5 fluid components viz., density, x-momentum, y-momentum,
z-momentum and energy. The last dimension 5 stores the modes of each variable. The first
mode corresponds to the zone centered value; these are the primal variables that we seek to
reconstruct, predict and evolve in time. The second, third and fourth modes correspond to the
piecewise linear variations in the x-, y- and z-directions.  These are obtained via a WENO limiter
step which provides the spatial reconstruction. The fifth mode corresponds to the temporal
dependence of the hyperbolic PDE and is provided by the ADER predictor step. Along with this,
there is another variable, which is a skinny version of the {\tt U} variable, called {\tt U\_skinny}, holds
the zone-centered primal variables of {\tt U}, and this gets updated for every timestep. The
{\tt U\_skinny} variable holds the minimal amount of data that is moved from the CPU (host) to
GPU (device) at the beginning of the time step and then back from the GPU (device) to CPU
(host) at the end of the timestep.

      The face-centered fluxes are stored in the variables {\tt Flux\_X}, {\tt Flux\_Y} and {\tt Flux\_Z}
respectively. For simplicity, they are all dimensioned as {\tt (0:nx, 0:ny, 0:nz)}, among which the x-
face information is stored in locations {\tt (0:nx, 1:ny, 1:nz)}, y-face information is stored in
locations {\tt (1:nx. 0:ny, 1:nz)} and the z-face information is stored in locations {\tt (1:nx. 1:ny, 0:nz)}.
For a zone located at {\tt i, j, k} its left and right x-faces are indexed as {\tt i-1, j, k} and {\tt i, j, k} 
respectively and the corresponding fluxes are stored in the {\tt Flux\_X} variable. In the same
way, its y-faces on the either are indexed as {\tt i, j-1, k} and {\tt i, j, k} and the corresponding fluxes
are stored in the {\tt Flux\_Y} variable. Similarly, its z-faces on the either are indexed as {\tt i, j, k-1}
and {\tt i, j, k} and the corresponding fluxes are stored in the {\tt Flux\_Z} variable. These fluxes will
be filled in via a Riemann solver-based corrector step.

      The time-update computed from these fluxes for each zone are stored in {\tt dU\_dt(nx, ny, nz, 5)}; 
its first three dimensions correspond to the x, y and z directions and the last dimension
5 holds the zone centered time-update corresponding to the 5 fluid variables. We also want to
create storage on the device for several useful scalars. These scalars include the current timestep
and the next timestep {\tt dt}, {\tt dt\_next}; the Courant Friedrichs Levy number, {\tt cfl} and the zone
sizes {\tt dx, dy, dz}. 
The directive {\tt !\$acc data create} allocates the memory for these variables on
the device. This is used to create the memory on the GPU for the variables {\tt U\_skinny}, {\tt U}, {\tt dU\_dt},
{\tt Flux\_X}, {\tt Flux\_Y}, {\tt Flux\_Z}. This directive does not copy the contents from CPU to GPU, it only
creates the memory space for these variables on the GPU.

      Once we have explained the data declaration and the use of various data variables in Fig.
1, it is valuable to describe the functional algorithmic steps in Fig. 1. It is expected that the
reader has a nicely modularized CFD code; with the result that we only point attention to the
additional steps that are needed to ready it for GPUs using OpenACC directives. Our goal, of
course, is to show how easy and transparent this upgrade is and to identify the minimal set of
OpenACC constructs needed to make this upgrade. In other words, we wish to provide the
techniques and tricks so that anyone with a nicely modular code can make the transition to
GPU-usage in the quickest and most painless way possible. The general calling sequence of the
subroutines involved in a HOG scheme can be seen inside the timestep loop, which starts from
{\tt DO i\_timestep = 1, number\_of\_timesteps} in Fig. 1. In this figure, the colorized text shows the
newly added codes which are necessary for the GPU implementation. The red color shows the
OpenACC directives, the blue color corresponds to the new data variable {\tt U\_skinny} and the
same blue shows the new subroutine(s) or code needed for the GPU implementation. The sequence
of the operations inside the time-loop are as follows:-

\begin{enumerate}[i.]
    \item Update the device with the zone-centered values stored in the {\tt U\_skinny} variable. This is
done with the {\tt !\$acc update device} directive in Fig. 1. We also copy some very useful scalars
from host to device within the same directive. This ensures that the absolute minimum amount of
data is moved from host to device. We will show later, in the results section, that not having this
minimalist data transfer strategy can result in a significant performance penalty.

\item Now that the minimal data has been moved over to the GPU, transfer the zone-centered
variables from {\tt U\_skinny} to the first mode of {\tt U}. This is done in {\tt U\_Skinny\_To\_U}.

\item Apply the necessary boundary conditions on {\tt U} via a call to {\tt Boundary\_Conditions}. This
ensures that physical values have been provided in the ghost zones.

\item reconstruct the x, y, z modes, i.e. {\tt U( :, :, :, :, 2:4)}, using a limiter. This is accomplished
using a call to Limiter. A nicely modularized HOG code should already have such a
subroutine/function so that we only need to show the minimal OpenACC upgrades needed in
such a limiter subroutine. (Indeed, the same upgrades would apply to a boundary condition
subroutine, which is why we dont repeat such information for a boundary condition subroutine.)

\item Using the spatial variation of the primal variables, predict the in-the-small temporal evolution
via a call to Predictor. This mode is returned by the predictor step as data that is stored in 
{\tt U( :, :, :, :, 5)}. Any modularly-written ADER-WENO code should have such a subroutine/function
so that we only need to show the minimal OpenACC upgrades needed in such a predictor
subroutine.

\item Now that the spatial and temporal modes are built, compute the face-centered HLL/HLLEM
fluxes via calls to Riemann solvers applied at the x-, y- and z-faces. This is done via calls to the
{\tt Make\_Flux\_X}, {\tt Make\_Flux\_Y}, {\tt Make\_Flux\_Z} subroutines. As before, any modularly-written
ADER-WENO code should have such a subroutine/function so that we only need to show the
minimal OpenACC upgrades needed in such a subroutine for evaluating fluxes.

\item Next, compute the time rate of change of the primal variables from the fluxes. Again, we
only need to show the minimal OpenACC upgrades that are needed.

\item Lastly, update the {\tt U( :, :, :, :, 1)} variable with {\tt dU\_dt}. On a CPU, this would have
completed the timestep. For a GPU implementation, we need to transfer the zone-centered values
back to the {\tt U\_skinny} variable which is to be used for transferring the minimum amount of data
back from device to host. This is done via the {\tt !\$acc update host} directive in Fig. 1. Using the
same directive, we also copy over the next prognosticated timestep, {\tt dt\_next}, for this particular
patch of data. In a multi-processor setting, this timestep will be minimized across different
patches on different nodes using an MPI reduction call; but providing that level of MPI-based
detail is not within the scope of this paper. Section V does, however, broadly discuss a scalability
study on a GPU-rich supercomputer where MPI was used for the communication across nodes.
\end{enumerate}

    In the following subsections, the small number of OpenACC changes that are needed for
each of the subroutines are discussed.

\subsubsection{Subroutine {\tt U\_Skinny\_to\_U( )}}

\begin{figure}
\begin{lstlisting}[language=Fortran]
SUBROUTINE U_Skinny_to_U( nx, ny, nz, U_skinny, U )

INTEGER nx, ny, nz
REAL U_skinny( -1:nx+2, -1:ny+2, -1:nz+2, 5)
REAL U( -1:nx+2, -1:ny+2, -1:nz+2, 5, 5)
INTEGER i, j, k

!$acc parallel vector_length(64) present( U_skinny, U )
!$acc loop gang vector collapse(3) independent
!$acc& private ( i, j, k )
DO k = -1, nz+2
  DO j = -1, ny+2
   DO i = -1, nx+2
     U( i, j, k, :, 1 ) = U_skinny( i, j, k, : )
   END DO
  END DO
END DO
!$acc end parallel

END SUBROUTINE U_Skinny_to_U
\end{lstlisting}
    \caption{shows the structure of the {\tt U\_Skinny\_to\_U} subroutine. 
    The array {\tt U\_skinny} contains the minimum data that should be carried 
    from the host to the device. This subroutine just copies the contents of 
    {\tt U\_skinny} to {\tt U( :, :, :, :, 1)}. To exploit the maximum available 
    amount of parallelism, it is imperative to collapse the three-fold nested loop.}
\end{figure}

    The {\tt U\_skinny} array stores the minimum amount of zone centered data, which is to be
transferred from the host to device and vice-versa for each time-step. This transfer of data is
essential in a parallel setting to provide the ghost zone information to the neighboring patches;
however, this MPI-based parallelism is not within the scope of this paper. Once the {\tt U\_skinny}
values are updated on the device using, {\tt !\$acc update device(U\_skinny)} at the beginning of
the time-step, the values are then copied to the first mode of the conserved variable {\tt U} in the
subroutine {\tt U\_skinny\_to\_U( )}, as shown in Fig. 2. The triply nested loop is for iterating over
the zones located along the x-, y-, and z-directions using the loop index variables {\tt i, j, k}. The
fourth index corresponds to the five fluid variables.
    There are mainly two OpenACC directives used in parallelizing the loops. The first one
declares the parallel region of the OpenACC. The beginning of the parallel region is marked with
{\tt !\$acc parallel} and the end of the parallel region is marked with {\tt !\$acc end parallel}. The
{\tt vector\_length(64)} specifies the length to be used in the vector operations on the GPU. Based
on our experience and testing with different GPUs a vector length 64 is chosen. (Depending on
the structure of their code, the reader should use a profiling tool to identify the optimal vector
length for their particular code.) The present clause is used to inform the GPU that the
variables that are included in the present clause are already declared or copied to the GPU.
Generally, the variables in the present clause are big-arrays which are shared between the loops
or smaller variables which remain constant during the loop execution. In this subroutine, the big-
arrays {\tt U\_skinny}, {\tt U} are placed inside the present clause.
    The second OpenACC directive, that is used here declares the loop and specifies the
nature of the parallelization used in the ensuing loop nest. Here, in order to utilize the
parallelization to the maximum extent all the 3-loops are collapsed and vectorized with the help
of the compiler. This is specified by the directive {\tt !\$acc loop gang vector collapse(3)}
independent. The last clause independent explicitly specifies to the compiler that the triply
nested loop iterations are data-independent and can be executed in parallel. In addition to this,
the small variables that are private to the loops can be explicitly specified here, using the
private clause. In this fashion, all of the shared and the private variables in a parallelizing loop
are explicitly specified to the compiler.
    Lastly, the end of the parallel region (which is also the end of the triply-nested loops) is
marked with the {\tt !\$acc end parallel}. In the next subsection directives used in the Limiter are
described.

\subsubsection{Subroutine {\tt Limiter( )}}

\begin{figure}
\begin{lstlisting}[language=Fortran]
SUBROUTINE Limiter( nx, ny, nz, U )
INTEGER nx, ny, nz
REAL U( -1:nx+2, -1:ny+2, -1:nz+2, 5, 5)
INTEGER i, j, k, n
REAL s1, s2, compression_factor, MC_Limiter
MC_Limiter( a, b, cfac ) = MIN(0.5 * ABS( a + b ), cfac * ABS( a ),
                              cfac * ABS( b )) * (SIGN(0.5, a ) + SIGN(0.5, b ))

!$acc data copyin( n, compression_factor )
DO n = 1, 5

  IF( n == 1 ) compression_factor = 2.0
  IF( n >  1 ) compression_factor = 1.5
!$acc update device( n, compression_factor )

!$acc parallel vector_length(64) present( n, U, compression_factor )
!$acc loop gang vector collapse(3) independent private( i, j, k, s1, s2 )
  DO k = 0, nz+1
    DO j = 0, ny+1
      DO i = 0, nx+1
        s1 = U( i+1, j, k, n, 1 ) - U(   i, j, k, n, 1 )
        s2 = U(   i, j, k, n, 1 ) - U( i-1, j, k, n, 1 )
        U( i, j, k, n, 2 ) = MC_Limiter( s1, s2, compression_factor )
      END DO
    END DO
  END DO
!$acc end parallel
END DO
! Do similarly for other modes in the y- and z-directions

!$acc end data 
END SUBROUTINE Limiter 
\end{lstlisting}
    \caption{shows the application of a simple Monotone-Centered limiter in the
    x-direction to all the fluid variables. The inline function {\tt MC\_Limiter} 
    is made to accommodate a compression factor which can indeed be changed 
    based on the flow field being limited. (For example, one might want to 
    provide a larger compression factor for the limiting of a density variable 
    than for the other variables. The limited x-slopes are stored in {\tt U( :, :, :, n, 2 )}.
    The limiter can also be applied in the y- and z-directions and the limited y- and z-slopes 
    are stored in {\tt U( :, :, :, n, 3)} and {\tt U( :, :, :, n, 4 )} .} 
\end{figure}

    Fig. 3 shows a simple monotone-centered limiter for the x-direction. For this an inline
function MC\_limiter is used as an example. It takes the slopes from either side and the
compression factor and gives back the limited x-slope. Apart from the triply nested DO-loop,
there is an outer-loop with a loop induction variable n, for each component of the fluid
dynamics. It can be noted that the value of n changes each time, the {\tt !\$acc parallel} region is
invoked. Also, the value of the {\tt compression\_factor} changes depending upon the fluid variable
to which the limiter is applied, i.e. n. (For example, one might want to use a steeper profile, i.e.
a higher compression factor, for the density variable so that the code can represent contact
discontinuities with reduced dissipation.) The efficient way to do this is to move the variables {\tt n,
compression\_factor} to the GPU and then update the device when there is a change in their
values. The first task is achieved by {\tt !\$acc data copyin}; this creates the memory for these
variables on the GPU and copies the values. Inside the outer loop given by {\tt DO n = 1,5}, once
the values of {\tt n, compression\_factor} are changed, they are then updated on the device using
{\tt !\$acc update device(n, compression\_factor)}. Once the use of these variables are complete,
they need to be deleted from the device using {\tt !\$acc end data}; this is done at the end of the
subroutine. Just to avoid any confusion about Fig. 3, please look at the last index in {\tt U}. It is the
second index that is updated using the first index. In other words, {\tt U(i, j, k, n, 2)}, which holds
the x-slope, is updated using {\tt U(i, j, k, n, 1)} and its neighbors, which indeed hold the primal
variables.
    The OpenACC directives for the inner 3-loops are similar to the previous description.
There are some specific points here which are worth noticing. In addition to the shared array
variable {\tt U}, the variables {\tt n, compression\_factor} are placed in the present clause. These two
variables remain constant for the parallel region. The private clause has two more real variables
{\tt a, b} which are private to the loop, and they are used in computing the x-slope. Each iteration of
the parallel loop nest shown in Fig. 3 will, therefore, have its own copy of the temporary
variables {\tt a, b} making them suitable for the computation of the undivided difference in that
figure. A similar privatization should be used for all variables that are used for the WENO
reconstruction at all other orders. Similarly, structured loops can be written for the y-, z-slopes of
the {\tt U} variable, completing the execution of second order accurate limiter. In the next
subsection, parallelization of the Predictor-subroutine is described.
 
\subsubsection{Subroutine {\tt Predictor( )}}
%
\begin{figure}
\begin{lstlisting}[language=Fortran]
SUBROUTINE Predictor( nx, ny, nz, U )
INTEGER nx, ny, nz
REAL U( -1:nx+2, -1:ny+2, -1:nz+2, 5, 5 )
INTEGER i, j, k
REAL U_ptwise(5, 5)    ! local variable for each zone

!$acc routine( Predictor_Ptwise ) seq

!$acc parallel vector_length(64) present( U )
!$acc loop gang vector collapse(3) independent
!$acc& private( i, j, k, U_ptwise )
DO k = 0, nz+1
  DO j = 0, ny+1
    DO i = 0, nx+1
      U_ptwise( 1:5, 1:4 ) = U( i, j, k, 1:5, 1:4)
      CALL Predictor_Ptwise( U_ptwise )            ! Obtain the time-mode
      U( i, j, k, 1:5, 5 ) = U_ptwise( 1:5, 5 )
    END DO
  END DO
END DO
!$acc end parallel
END SUBROUTINE Predictor

SUBROUTINE Predictor_Ptwise( U_ptwise )
!$acc routine seq
! Pointwise predictor step
END SUBROUTINE Predictor_Ptwise
\end{lstlisting}
    \caption{shows the predictor step for building the 5th mode, i.e., the mode 
    that contains the time rate of change. The zone centered value {\tt U( :, :, :, 1:5, 1 )} 
    and the space-modes {\tt U( :, :, :, 1:5, 2:4 )} that are obtained from the limiter are 
    copied to a zone-wise variable called {\tt U\_ptwise} for all the 5-components 
    of the hydrodynamics. Using these and the expression of fluxes the 
    {\tt Predictor\_Ptwise} would compute the 5th mode, which contains the temporal 
    evolution. This is then copied back to the 5th mode of the variable 
    {\tt U( :, :, :, 1:5, 5 )} for all the 5 components of the hydrodynamics.
    Here, in the {\tt !\$acc routine( \dots$~$) seq}, the OpenACC declaration of the 
    subroutine can be noted. It declares to the GPU that the subroutine would 
    be called in a sequential manner. A similar statement needs to be made 
    inside the {\tt Predictor\_Ptwise} subroutine and this is shown in the bottom 
    part of the figure.}
\end{figure}

    The predictor subroutine builds the fifth mode that contains the in-the-small time rate of
change in the solution vector for the PDE. The general structure of the predictor subroutine,
along with OpenACC directives, is shown in Fig. 4. For the higher order schemes, the predictor
step can be quite long and complicated; therefore, it is commonly packaged into a subroutine,
which would operate once on each zone during each timestep. (We refer to it as a zone-wise
subroutine.) The zone-wise subroutine is then called inside the triply nested loop. Its inputs are
the first four modes of the five fluid variables, in other words, the zone-wise input 
{\tt U\_ptwise(1:5, 1:4) = U( i, j, k, 1:5, 1:4)}. 
Using the zone centered value and the spatial modes the {\tt Predictor\_Ptwise} subroutine 
builds the fifth mode; i.e., the time rate of change for all five
fluid variables. This is then stored into the fifth mode of {\tt U}, i.e. in 
{\tt U( i, j, k,1:5, 5)}.
    The OpenACC directives around the loop are exactly similar to that of the previous
subroutines. It can be noted that the variable {\tt U\_ptwise} is a small 5x5 array, and it is
classified as the private variable. (Incidentally, the use of such very small arrays also makes the
pointwise predictor very cache friendly towards the small caches on the GPU.) By default, the
OpenACC classifies an array variable as a shared variable, which is not true here. It is imperative
that the {\tt U\_ptwise} is explicitly identified as a private variable for the GPU. The subroutine
call inside the triply nested loop needs to be specified to the GPU as a sequential call, without
any parallelization. This is accomplished with the help of the {\tt !\$acc routine seq} directive. This
directive is placed inside the {\tt Predictor\_Ptwise} subroutine, as shown in the bottom of Fig 4.
In addition to this, the subroutines may be placed in different files during the compilation,
therefore an additional mention of this sequential calling needs to be present inside the
Predictor subroutine. This is accomplished by the OpenACC directive in the third line of the
Fig. 4, that is, {\tt !\$acc routine(Predictor\_Ptwise) seq}. It informs the compiler that the
{\tt Predictor\_Ptwise} subroutine is a purely sequential subroutine. In the next subsection, the
flux evaluation subroutine is discussed.

\subsubsection{Subroutine {\tt Make\_Flux\_X( )}}
\begin{figure}
\begin{lstlisting}[language=Fortran]
SUBROUTINE Make_Flux_X( nx, ny, nz, U, Flux_X )
INTEGER nx, ny, nz
REAL U( -1:nx +2, -1:ny+2, -1:nz+2, 5, 5 )
REAL Flux_X( 0:nx, 0:ny, 0:nz, 5 )
INTEGER i, j, k
REAL flux_x_ptwise(5), U_L( 5), U_R(5)

!$acc routine( Riemann_Solver_Ptwise ) seq

!$acc parallel vector_length(64) present( U, Flux_X )
!$acc loop gang vector collapse (3) independent
!$acc&      private( i, j, k, U_L, U_R, flux_x_ptwise )
DO k = 1, nz
  DO j = 1, ny
    DO i = 0, nx
      U_L(:) = U( i, j, k, :, 1 ) + 0.5 * U( i, j, k, :, 2 )
                                  + 0.5 * U( i, j, k, :, 5 )

      U_R(:) = U( i+1, j, k, :, 1 ) - 0.5 * U( i+1, j, k, :, 2 )
                                    + 0.5 * U( i+1, j, k, :, 5 )

      CALL Riemann_Solver_Ptwise( U_ L, U_ R, flux_x_ptwise )
      Flux_X( i, j, k, :) = flux_x_ptwise(:)
    END DO
  END DO
END DO
!$acc end parallel
END SUBROUTINE Make_Flux_X

SUBROUTINE Riemann_Solver_Ptwise( U_L, U_R, flux_x_ptwise )
!$acc routine seq
! Build the HLL/HLLEM flux at the x -face
END SUBROUTINE Riemann_Solver_Ptwise
\end{lstlisting}
    \caption{shows the pseudocode for the subroutine to build the flux along the 
    x-faces. The {\tt U\_L} variable stores the value of {\tt U} obtained from the zone which
    is at the left to the face and the {\tt U\_R} variable stores the value of {\tt U} obtained
    from the zone which is at the right to the face.  Thus the facial values of the
    conserved variable {\tt U} are stored in {\tt U\_L, U\_R}. This is then fed in to the Riemann
    solver subroutine to build the flux for the corresponding face. The Riemann solver
    subroutine is called in a serial fashion. The same has been indicated to the
    GPU in the 3rd line of the figure. It declares to the GPU that the subroutine
    would be called in a sequential manner. A similar statement need to be made inside
    the {\tt Riemann\_Solver\_Ptwise} subroutine and this is shown in the bottom part of
    the figure. The pointwise fluxes obtained from the Reimann solver are stored
    in {\tt Flux\_X ( 0:nx, 1:ny, 1:nz, :)}. 
    The subroutines for building the {\tt Flux\_Y ( 1:nx, 0:ny, 1:nz, :), 
    Flux\_Z ( 1:nx, 1:ny, 0:nz, :)} can be written in the same way, as that of the
    subroutine shown here. For the y- and z-fluxes, instead of the 2nd mode, the
    3rd and 4th modes of {\tt U} will be used for building the ‘y’ and ‘z’ face fluxes respectively.}
\end{figure}

    The previous predictor step builds the spatial and temporal modes. Using those modes,
the face-centered HLL/HLLEM fluxes can be built by calling the Riemann solvers at the x-, y-
and z-faces. This is done via calls to the {\tt Make\_Flux\_X}, {\tt Make\_Flux\_Y}, {\tt Make\_Flux\_Z}
subroutines. As an example, the general structure of the {\tt Make\_Flux\_X} subroutine is shown in
Fig. 5. Inside the triply nested loop, the conserved variables from the left ({\tt U\_L}) and right ({\tt U\_R})
zones at the face are constructed and sent to the {\tt Riemann\_Solver\_Ptwise} to compute the
resolved flux at the face. The values of the resolved flux are then copied to the bigger array
{\tt Flux\_X} from the small pointwise array {\tt flux\_x\_ptwise} for all the five fluid dynamic
variables. (As with the predictor step, this pointwise processing ensures that each call to
{\tt Riemann\_Solver\_Ptwise} is independent of any other call. It also ensures that
{\tt Riemann\_Solver\_Ptwise} only uses very small private arrays making it very cache friendly
towards the small caches on the GPU.)
    As mentioned before, by default, the OpenACC classifies all array variables as shared
variables, which is not the case for {\tt U\_L}, {\tt U\_R} and {\tt flux\_x\_ptwise}. Therefore, it is
imperative that the {\tt U\_L}, {\tt U\_R} and {\tt flux\_x\_ptwise} are explicitly identified as private
variables for the GPU. The same can be seen above the triply nested loop where all of the above-
mentioned variables, as well as the indices {\tt i, j, k}, are declared as private to the compiler. The
pointwise variables {\tt U\_L}, {\tt U\_R} are sent as input to the Riemann solver and the
{\tt flux\_x\_ptwise} is the output from the Riemann solver. The {\tt Riemann\_Solver\_Ptwise}
subroutine is executed sequentially inside the triply nested loop. Hence the OpenACC directive
{\tt !\$acc routine seq} is placed inside the subroutine as shown in the bottom of Fig. 5. In addition
to this the same is conveyed to the {\tt Make\_Flux\_X} subroutine in the third line; using the
directive {\tt !\$acc routine(Riemann\_Solver\_Ptwise) seq}.  In the next subsection the building of
the time rate of change of the conserved variable {\tt U} is discussed.

\subsubsection{Subroutine {\tt Make\_dU\_dt( )}}

\begin{figure}
\begin{lstlisting}[language=Fortran]
SUBROUTINE Make_dU_dt( nx, ny, nz, Flux_X, Flux_Y, Flux_Z,
                       dU_dt, dt, dx, dy, dz )
INTEGER nx, ny, nz
REAL U( -1:nx+2, -1:ny+2, -1:nz+2, 5, 5)
REAL, DIMENSION( 0:nx, 0:ny, 0:nz, 5 ) :: Flux_X, Flux_Y, Flux_Z
REAL dU_dt( nx, ny, nz, 5 )
REAL dt, dx, dy, dz
INTEGER i, j, k

!$acc parallel vector_length(64) present (
!$acc& dt, dx, dy, dz, Flux_X, Flux_Y, Flux_Z, dU_dt )
!$acc loop gang vector collapse(3) independent private( i, j, k )
DO k = 1, nz
  DO j = 1, ny
    DO i = 1, nx
      dU_dt( i, j, k, :) =
       - dt * (Flux_X( i, j, k, : ) - Flux_X( i-1, j, k, :)) / dx
       - dt * (Flux_Y( i, j, k, : ) - Flux_Y( i, j-1, k, :)) / dy
       - dt * (Flux_Z( i, j, k, : ) - Flux_Z( i, j, k-1, :)) / dz
    END DO
  END DO
END DO
!$acc end parallel

END SUBROUTINE Make_dU_dt
\end{lstlisting}
    \caption{hows the subroutine to build the dU\_dt update from the fluxes obtained
    from the {\tt Make\_Flux\_X, Make\_Flux\_Y and Make\_Flux\_Z} subroutines. The 
    {\tt dU\_dt ( 1:nx, 1:ny, 1:nz, :)} variable is computed for all the active zones
    of the simulation. The common variables are present on the GPU and they are indicated
    as {\tt present( \dots$~$)}; the local variables are mentioned in the 
    {\tt private( \dots$~$)} part.}
\end{figure}

    The next step would be to build the time rate of change for the conserved variable U
from the facial fluxes. This is carried out in the {\tt Make\_dU\_dt} subroutine, as shown in Fig. 6. In
addition to the facial flux variables, {\tt Flux\_X}, {\tt Flux\_Y} and {\tt Flux\_Z}, the zone sizes {\tt dx, dy, dz}
and the timestep size {\tt dt} are needed to compute the evolution of {\tt U}. The variables {\tt Flux\_X},
{\tt Flux\_Y}, {\tt Flux\_Z}, {\tt dU\_dt} are already present on the GPU memory. The additional variables, 
{\tt dt, dx, dy, dz} are now needed on the GPU; 
however, realize that they were already copied from the
host to device, using the directive {\tt !\$acc update device( \dots, dt, cfl, dx, dy, dz)} in Fig. 1. Once
moved to the device, all of these variables can be placed inside the present clause as shared
variables for the GPU. In the next subsection, the subroutine which updates the {\tt U} with the
time rate of change {\tt dU\_dt} is discussed.

\subsubsection{Subroutine {\tt Update\_U\_Timestep( )}}
\begin{figure}
\begin{lstlisting}[language=Fortran]
SUBROUTINE Update_U_Timestep ( nx, ny, nz, U, U_skinny, dU_dt,
                               dt_next, cfl, dx, dy, dz )
INTEGER nx, ny, nz
REAL U( -1:nx+2, -1:ny+2, -1:nz+2, 5, 5)
REAL U_skinny( -1:nx+2, -1:ny+2, -1:nz+2, 5 ), dU_dt( nx, ny, nz, 5 )
REAL dt_next, cfl, dx, dy, dz
INTEGER i, j, k
REAL U_ptwise(5), dt1

!$acc routine( Eval_Tstep_Ptwise ) seq

dt_next = 1.0e32
!$acc update device( dt_next )

!$acc parallel vector_length(64) present ( U,
!$acc&      U_skinny, dU_dt, cfl, dx, dy, dz )
!$acc loop gang vector collapse(3) independent
!$acc&     private( i, j, k, U_ptwise, dt1 ) 
!$acc&     REDUCTION( MIN:dt_next )
DO k = 1, nz
  DO j = 1, ny
    DO i = 1, nx

      U( i, j, k, :, 1 ) = U( i, j, k, :, 1 ) + dU_dt( i, j, k, : )
      U_skinny( i, j, k, : ) = U( i, j, k, :, 1 )
      U_ptwise(:) = U( i, j, k, :, 1 )

      CALL Eval_Tstep_Ptwise( cfl, dx, dy, dz, U_ptwise, dt1 )
      dt_next = MIN(dt_next, dt1 )

    END DO
  END DO
END DO
!$acc end parallel

END SUBROUTINE Update_U_Timestep
\end{lstlisting}
    \caption{shows the subroutine to update the {\tt U} from the {\tt dU\_dt}. In addition
    to this, here the {\tt U\_skinny} is updated from the new zone-centered values of
    {\tt U}. The minimum amount of data stored in the {\tt U\_skinny} will be moved back to
    the CPU, demonstrating the importance of minimal data movement in a CPU+GPU computer.
    In a parallel setting, this data can be used to provide ghost zone information
    to neighboring patches; but that is not within the scope of this paper.
    Also in this function, the next-timestep {\tt dt\_next} is estimated, using the CFL
    number, zone size and U. The smallest among all the zones is found using the
    {\tt MIN()} function. In order to correctly parallelize this reduction operation
    the OpenACC command, {\tt REDUCTION( MIN: dt\_next )} is used.}
\end{figure}

      In this last step, the {\tt U( :, :, :, :, 1)} variable is updated with {\tt dU\_dt} and the next
prognosticated timestep, {\tt dt\_next} is evaluated. The general structure of this subroutine is shown
in Fig. 7. The triply nested loop is used as before and the five fluid dynamic variables are
updated for each zone. Notice that {\tt cfl, dx, dy, dz} (which are all needed for evaluating the next
timestep via {\tt Eval\_Tstep\_Ptwise}) are already present on the device. Since {\tt dt\_next} is
reinitialized in this subroutine, it is updated on the device. The {\tt REDUCTION( MIN: dt\_next )}
pragma tells the compiler to treat {\tt dt\_next} as a reduction variable that has to be reduced over all
the cores of the GPU. This completes our description of the timestep.

      As mentioned before, there is a need for the transfer of information from the device to
host at each time-step. This is necessary during a parallel execution across many nodes, where
the ghost cell information must be communicated (using MPI-based messaging) from one patch
to the other at each time-step. Even though, this is not within the scope of this paper, for the sake
of completeness (and as a demonstration of a very popular use case) the minimum amount of
data is packed into {\tt U\_skinny} and moved from host to device via the 
{\tt !\$acc update host(U\_skinny)} directive in Fig. 1. To prepare for this, 
the updated zone centered values of {\tt U} are
copied to the {\tt U\_skinny} variable in this subroutine, inside the triply nested loop. The
prognosticated value for the next timestep, which should eventually be minimized across nodes
in a parallel implementation, is also ready to be sent over from the device to the host.

\section{Results}
    In this section we present results associated with GPU speedup for the CFD, MHD and
CED codes. All three codes are based on higher order Godunov scheme philosophy. In each
case, we compare performance on an NVIDIA A100 PCIe GPU with 80GB of memory relative
to the performance of a single core of a 3.0GHz Xeon Gold 6248R CPU. We used the NVIDIA
NVHPC (version 22.2) Fortran compiler for all our runs with an optimization level of 3. This
compiler works on CPUs and it also has full-featured OpenACC which enables it to work on
GPUs.

\subsection{CFD Results}
\begin{table}[!ht]
    \caption{A100 GPU speedup for the CFD code when the {\tt U\_skinny} trick was not used.}
    \vspace{2mm}
    \centering
    \begin{tabular}{|c|c|c|c|c|c|c|c|}
    \hline
              &                 & Zones /         &              & Zones /      &          \\ 
        Zones & Serial time (s) & Serial time (s) & GPU time (s) & GPU time (s) & Speed up \\ \hline
        O2 & ~ & ~ & ~ & ~ & ~ \\ \hline
        $48^3$ & 1.43E-01 & 7.75E+05 & 8.43E-03 & 1.31E+07 & 16.9 \\ \hline
        $96^3$ & 1.14E+00 & 7.75E+05 & 3.98E-02 & 2.22E+07 & 28.7 \\ \hline
        $144^3$ & 3.67E+00 & 8.14E+05 & 1.10E-01 & 2.72E+07 & 33.5 \\ \hline
        $192^3$ & 8.73E+00 & 8.11E+05 & 2.62E-01 & 2.70E+07 & 33.3 \\ \hline
        ~ & ~ & ~ & ~ & ~ & ~ \\ \hline
        O3 & ~ & ~ & ~ & ~ & ~ \\ \hline
        $48^3$ & 3.37E-01 & 3.29E+05 & 2.56E-02 & 4.32E+06 & 13.1 \\ \hline
        $96^3$ & 2.41E+00 & 3.67E+05 & 1.08E-01 & 8.18E+06 & 22.3 \\ \hline
        $144^3$ & 7.85E+00 & 3.81E+05 & 3.31E-01 & 9.03E+06 & 23.7 \\ \hline
        $192^3$ & 1.86E+01 & 3.81E+05 & 7.99E-01 & 8.86E+06 & 23.2 \\ \hline
        ~ & ~ & ~ & ~ & ~ & ~ \\ \hline
        O4 & ~ & ~ & ~ & ~ & ~ \\ \hline
        $48^3$ & 1.13E+00 & 9.75E+04 & 5.75E-02 & 1.92E+06 & 19.7 \\ \hline
        $96^3$ & 7.04E+00 & 1.26E+05 & 3.43E-01 & 2.58E+06 & 20.5 \\ \hline
        $144^3$ & 2.23E+01 & 1.34E+05 & 1.14E+00 & 2.61E+06 & 19.5 \\ \hline
        $192^3$ & 5.14E+01 & 1.38E+05 & 2.97E+00 & 2.39E+06 & 17.3 \\ \hline
    \end{tabular}
\end{table}

\begin{table}[!ht]
    \caption{A100 GPU speedup for the CFD code when the {\tt U\_skinny} trick was used.}
    \vspace{2mm}
    \centering
    \begin{tabular}{|c|c|c|c|c|c|c|c|}
    \hline
              &                 & Zones /         &              & Zones /      &          \\ 
        Zones & Serial time (s) & Serial time (s) & GPU time (s) & GPU time (s) & Speed up \\ \hline
        O2 & ~ & ~ & ~ & ~ & ~ \\ \hline
        $48^3$ & 1.43E-01 & 7.73E+05 & 3.43E-03 & 3.22E+07 & 41.7 \\ \hline
        $96^3$ & 1.13E+00 & 7.86E+05 & 1.58E-02 & 5.60E+07 & 71.3 \\ \hline
        $144^3$ & 3.71E+00 & 8.04E+05 & 4.62E-02 & 6.46E+07 & 80.4 \\ \hline
        $192^3$ & 9.11E+00 & 7.77E+05 & 1.08E-01 & 6.55E+07 & 84.4 \\ \hline
        ~ & ~ & ~ & ~ & ~ & ~ \\ \hline
        O3 & ~ & ~ & ~ & ~ & ~ \\ \hline
        $48^3$ & 3.32E-01 & 3.33E+05 & 5.17E-03 & 2.14E+07 & 64.3 \\ \hline
        $96^3$ & 2.47E+00 & 3.58E+05 & 2.47E-02 & 3.58E+07 & 100 \\ \hline
        $144^3$ & 7.90E+00 & 3.78E+05 & 7.30E-02 & 4.09E+07 & 108.2 \\ \hline
        $192^3$ & 1.88E+01 & 3.77E+05 & 1.72E-01 & 4.12E+07 & 109.3 \\ \hline
        ~ & ~ & ~ & ~ & ~ & ~ \\ \hline
        O4 & ~ & ~ & ~ & ~ & ~ \\ \hline
        $48^3$ & 1.11E+00 & 9.98E+04 & 1.70E-02 & 6.49E+06 & 65 \\ \hline
        $96^3$ & 7.14E+00 & 1.24E+05 & 9.24E-02 & 9.58E+06 & 77.4 \\ \hline
        $144^3$ & 2.24E+01 & 1.33E+05 & 2.85E-01 & 1.05E+07 & 78.7 \\ \hline
        $192^3$ & 5.21E+01 & 1.36E+05 & 6.21E-01 & 1.14E+07 & 83.8 \\ \hline
    \end{tabular}
\end{table}

    Tables I and II show the speedups associated with a CFD code that is run on an A100
GPU. In each case we compare the single core CPU performance with the GPU performance.
(Section V will show cases where the speeds of several GPUs are compared to the speeds
equivalent numbers of multicore CPUs for a real-world application.) Table I shows the results
when the {\tt U\_skinny} trick is not used; Table II shows the results when the {\tt U\_skinny} trick is
indeed used. Table I shows the traditional amount of speedup that is often reported by several
practitioners for GPUs; and indeed it is fair to say that it is encouraging but not very impressive.
We see a dramatic improvement in the results shown in Table II; and we see that the speedup is
indeed quite compelling. This shows the importance of only making minimal data motion
between the host and the device. From Table II we also see that the third order scheme shows
better speedup on comparable meshes than the second order scheme. This is an advantage that
we will see persisting even in MHD calculations. The CFD code and the MHD code use similar
type of WENO, ADER and Riemann solver algorithms, therefore, their speedup trend lines are
similar. Table II also demonstrates that with all the right algorithmic choices, along with deft
utilization of OpenACC pragmas, it is indeed possible to unlock the potential of GPU computing
for CFD applications.

    It is important to point out another issue. Please look at the raw numbers for zones
updated per second on the CPU in Table II. We see that at second order, the CFD code updates
~0.77 million zones per second on a CPU. At third order, the CFD code updates ~0.377 million
zones per second on a CPU. At fourth order, the CFD code updates ~0.136 million zones per
second on a CPU. These are very good speeds for a CFD code on a CPU, showing that the CPU
code was already highly optimized for CPUs. On CPUs, we have made efforts to ensure that the
code is optimally cache-friendly. It is harder to show good GPU speedups when the code is
already optimized for CPUs. The fact that Table II, nevertheless, shows good speedup on a GPU
means that we had to do everything right even on the GPU in order to extract that level of
speedup.

\begin{table}[!ht]
    \centering
    \caption{comparing the fractional time spent in the ADER subroutine at various 
    orders on CPU and GPU for the CFD code.}
    \vspace{2mm}
    \begin{tabular}{|c|c|c|c|c|}
    \hline
              &                 & Fraction of time         &              & Fraction of time         \\ 
        Zones & Serial time (s) &                  in ADER & GPU time (s) &                  in ADER \\ \hline
        O2 & ~ & ~ & ~ & ~ \\ \hline
        $48^3$ & 4.48E-02 & 0.313 & 3.49E-04 & 0.102 \\ \hline
        $96^3$ & 2.82E-01 & 0.25 & 1.99E-03 & 0.126 \\ \hline
        $144^3$ & 8.76E-01 & 0.236 & 6.13E-03 & 0.133 \\ \hline
        $192^3$ & 2.03E+00 & 0.222 & 1.40E-02 & 0.13 \\ \hline
        ~ & ~ & ~ & ~ & ~ \\ \hline
        O3 & ~ & ~ & ~ & ~ \\ \hline
        $48^3$ & 1.62E-01 & 0.488 & 1.76E-03 & 0.34 \\ \hline
        $96^3$ & 1.10E+00 & 0.447 & 1.07E-02 & 0.433 \\ \hline
        $144^3$ & 3.42E+00 & 0.433 & 3.30E-02 & 0.451 \\ \hline
        $192^3$ & 7.87E+00 & 0.419 & 7.52E-02 & 0.437 \\ \hline
        ~ & ~ & ~ & ~ & ~ \\ \hline
        O4 & ~ & ~ & ~ & ~ \\ \hline
        $48^3$ & 6.19E-01 & 0.559 & 1.29E-02 & 0.755 \\ \hline
        $96^3$ & 3.78E+00 & 0.528 & 7.32E-02 & 0.792 \\ \hline
        $144^3$ & 1.12E+01 & 0.501 & 2.20E-01 & 0.772 \\ \hline
        $192^3$ & 2.56E+01 & 0.491 & 4.91E-01 & 0.791 \\ \hline
    \end{tabular}
\end{table}

    Comparing the speedups at second and third order in Table II, we see that the third order
scheme shows substantially better speedup on the GPU than on the CPU. This is because the data
motion from CPU to GPU is the same for both schemes, however, the third order scheme uses
more floating point operations per zone than the second order scheme. The fourth order scheme
has a comparable speedup to the second order scheme, but the speedup is not so good when
compared to the third order scheme. This is explained in Table III where we focus exclusively on
the fraction of time that the CFD code spends within the ADER predictor step. At second and
third orders, the fraction of time spent within the ADER predictor step is comparable on the CPU
and on the GPU. Admittedly, the third order ADER takes the place of a three stage third order
Runge-Kutta scheme while the second order ADER scheme takes the place of a two stage second
order Runge-Kutta scheme. As a result, the fraction of time spent in the ADER step is larger for
the third order scheme compared to the second order scheme. But now please compare the
fractional times spent in the ADER predictor step at fourth order on CPU and GPU. We see that
the GPU code spends a much larger fraction of its time in the ADER predictor step compared to
the CPU code. The physical reason for that is because ADER predictor steps use increasingly
larger amounts of memory with increasing order. The CPU has a large cache, so it works fine
with higher order ADER. However, the GPU cores have a very small cache, so the fourth order
ADER scheme is not as adept at using the (smaller) cache on the GPU cores very efficiently. At
the point of this writing, we do not know of any strategies to reduce the cache usage of higher
order ADER predictor steps. Even so, this work serves to identify one of the bottlenecks in
present-day GPU computing.

\begin{table}[!ht]
    \centering
    \caption{A100 GPU speedup for the CFD code where we exclusively focus on 
    the performance differences between the ADER-based scheme and the Runge-Kutta-based scheme.}
    \vspace{2mm}
    \begin{tabular}{|c|c|c|c|c|c|c|}
    \hline
              &           &           &           &           & Relative  &    ~     \\
              &   ADER -  &           &           & Speed up  & Speed up  & Speed up \\
              &    CPU    & RK - CPU  & RK - GPU  & ADER-GPU/ &  RK-GPU/  &  RK-GPU/ \\
        Zones & Zones/sec & Zones/sec & Zones/sec & ADER-CPU  & ADER-CPU  &  RK-CPU  \\ \hline
        O2 & ~ & ~ & ~ & ~ & ~ & ~ \\ \hline
        $48^3$ & 7.73E+05 & 5.56E+05 & 1.61E+07 & 41.7 & 20.88 & 29.04 \\ \hline
        $96^3$ & 7.86E+05 & 5.02E+05 & 3.08E+07 & 71.3 & 39.15 & 61.29 \\ \hline
        $144^3$ & 8.04E+05 & 5.26E+05 & 3.60E+07 & 80.4 & 44.72 & 68.43 \\ \hline
        $192^3$ & 7.77E+05 & 4.99E+05 & 3.66E+07 & 84.4 & 47.08 & 73.33 \\ \hline
        ~ & ~ & ~ & ~ & ~ & ~ & ~ \\ \hline
        O3 & ~ & ~ & ~ & ~ & ~ & ~ \\ \hline
        $48^3$ & 3.33E+05 & 1.92E+05 & 1.03E+07 & 64.3 & 30.89 & 53.58 \\ \hline
        $96^3$ & 3.58E+05 & 2.01E+05 & 1.86E+07 & 100 & 51.83 & 92.58 \\ \hline
        $144^3$ & 3.78E+05 & 2.04E+05 & 2.08E+07 & 108.2 & 54.95 & 101.73 \\ \hline
        $192^3$ & 3.77E+05 & 1.99E+05 & 2.05E+07 & 109.3 & 54.42 & 102.87 \\ \hline
        ~ & ~ & ~ & ~ & ~ & ~ & ~ \\ \hline
        O4 & ~ & ~ & ~ & ~ & ~ & ~ \\ \hline
        $48^3$ & 9.98E+04 & 3.98E+04 & 5.05E+06 & 65 & 50.59 & 126.97 \\ \hline
        $96^3$ & 1.24E+05 & 4.68E+04 & 8.43E+06 & 77.4 & 68.04 & 179.95 \\ \hline
        $144^3$ & 1.33E+05 & 4.72E+04 & 9.67E+06 & 78.7 & 72.67 & 204.95 \\ \hline
        $192^3$ & 1.36E+05 & 4.75E+04 & 9.59E+06 & 83.8 & 70.51 & 201.65 \\ \hline
    \end{tabular}
\end{table}

    The curious reader might still want to know what would happen if we had used an
equivalent Runge-Kutta temporal update strategy. This is explored in Table IV. (Please recall
that a second order Runge-Kutta scheme uses two sub-steps; a third order Runge-Kutta scheme
uses three sub-steps but a fourth order Runge-Kutta scheme uses five sub-steps.) The fifth
column in Table IV just repeats the last column from Table II. The sixth column compares
Runge-Kutta speed on a GPU compared to ADER speed on a CPU. The seventh column in Table
IV compares the Runge-Kutta speed on a GPU to the Runge-Kutta speed on a CPU. Comparing
the fifth and sixth columns in Table IV, we see that the speedup of a Runge-Kutta based code on
the GPU relative to the CPU speed of an ADER-based code is not very impressive. The last
column of Table IV compares the speedup of a Runge-Kutta based code on the GPU relative to
the same code on a CPU. Then indeed we see quite wonderful speedups in the last column of
Table IV. There are two ways to interpret this column. First, if all that the user has is a Runge-
Kutta code, then this column shows that the GPU speeds will be vastly better than CPU speeds,
as long as the {\tt U\_skinny} trick is used. This is a strong positive and gives us a pathway to
dramatically improving the speeds of Runge-Kutta based codes. Second, we realize that the
Runge-Kutta timestepping algorithm does not use an ADER step, which was found to be a
bottleneck in Table III. Therefore, the last column of Table IV also shows us that with increasing
order, and all the way up to fourth order, a Runge-Kutta based code will show improving
speedup with increasing order. However, note that the intrinsic number of zones updated per
second by a Runge-Kutta based CFD code is consistently lower compared to an ADER based
CFD code. (Compare the fourth column in Table IV to the fifth column in Table II.) This trend is
true for CPUs and GPUs and shows us that a one-time investment in ADER algorithms can be
worthwhile.

\subsection{MHD Results}

\begin{table}[!ht]
    \caption{A100 GPU speedup for the MHD code when the {\tt U\_skinny} trick was not used.}
    \vspace{2mm}
    \centering
    \begin{tabular}{|c|c|c|c|c|c|c|c|}
    \hline
              &                 & Zones /         &              & Zones /      &          \\ 
        Zones & Serial time (s) & Serial time (s) & GPU time (s) & GPU time (s) & Speed up \\ \hline
        O2 & ~ & ~ & ~ & ~ & ~ \\ \hline
        $48^3$ & 8.23E-01 & 1.34E+05 & 2.53E-02 & 4.37E+06 & 32.5 \\ \hline
        $96^3$ & 6.27E+00 & 1.41E+05 & 1.07E-01 & 8.24E+06 & 58.4 \\ \hline
        $144^3$ & 2.08E+01 & 1.44E+05 & 3.35E-01 & 8.92E+06 & 62.1 \\ \hline
        $192^3$ & 5.12E+01 & 1.38E+05 & 7.36E-01 & 9.62E+06 & 69.6 \\ \hline
        ~ & ~ & ~ & ~ & ~ & ~ \\ \hline
        O3 & ~ & ~ & ~ & ~ & ~ \\ \hline
        $48^3$ & 1.35E+00 & 8.19E+04 & 5.09E-02 & 2.17E+06 & 26.5 \\ \hline
        $96^3$ & 1.08E+01 & 8.16E+04 & 2.76E-01 & 3.21E+06 & 39.4 \\ \hline
        $144^3$ & 3.72E+01 & 8.04E+04 & 7.85E-01 & 3.81E+06 & 47.3 \\ \hline
        $192^3$ & 9.25E+01 & 7.66E+04 & 1.69E+00 & 4.19E+06 & 54.7 \\ \hline
        ~ & ~ & ~ & ~ & ~ & ~ \\ \hline
        O4 & ~ & ~ & ~ & ~ & ~ \\ \hline
        $48^3$ & 3.06E+00 & 3.62E+04 & 1.35E-01 & 8.21E+05 & 22.7 \\ \hline
        $96^3$ & 2.28E+01 & 3.87E+04 & 7.34E-01 & 1.21E+06 & 31.1 \\ \hline
        $144^3$ & 7.86E+01 & 3.80E+04 & 2.57E+00 & 1.16E+06 & 30.6 \\ \hline
        $192^3$ & 1.89E+02 & 3.75E+04 & 4.25E+00 & 1.66E+06 & 44.3 \\ \hline
    \end{tabular}
\end{table}

\begin{table}[!ht]
    \caption{A100 GPU speedup for the MHD code when the {\tt U\_skinny} trick was used.}
    \vspace{2mm}
    \centering
    \begin{tabular}{|c|c|c|c|c|c|c|c|}
    \hline
              &                 & Zones /         &              & Zones /      &          \\ 
        Zones & Serial time (s) & Serial time (s) & GPU time (s) & GPU time (s) & Speed up \\ \hline
        O2 & ~ & ~ & ~ & ~ & ~ \\ \hline
        $48^3$ & 8.01E-01 & 1.38E+05 & 1.66E-02 & 6.65E+06 & 48.2 \\ \hline
        $96^3$ & 6.18E+00 & 1.43E+05 & 6.71E-02 & 1.32E+07 & 92.2 \\ \hline
        $144^3$ & 2.09E+01 & 1.43E+05 & 2.13E-01 & 1.40E+07 & 98.1 \\ \hline
        $192^3$ & 5.26E+01 & 1.34E+05 & 4.66E-01 & 1.52E+07 & 112.9 \\ \hline
        ~ & ~ & ~ & ~ & ~ & ~ \\ \hline
        O3 & ~ & ~ & ~ & ~ & ~ \\ \hline
        $48^3$ & 1.33E+00 & 8.29E+04 & 2.33E-02 & 4.74E+06 & 57.2 \\ \hline
        $96^3$ & 1.05E+01 & 8.41E+04 & 1.01E-01 & 8.74E+06 & 103.9 \\ \hline
        $144^3$ & 3.81E+01 & 7.83E+04 & 3.18E-01 & 9.38E+06 & 119.8 \\ \hline
        $192^3$ & 9.43E+01 & 7.50E+04 & 7.59E-01 & 9.33E+06 & 124.3 \\ \hline
        ~ & ~ & ~ & ~ & ~ & ~ \\ \hline
        O4 & ~ & ~ & ~ & ~ & ~ \\ \hline
        $48^3$ & 3.09E+00 & 3.58E+04 & 5.16E-02 & 2.15E+06 & 60 \\ \hline
        $96^3$ & 2.28E+01 & 3.87E+04 & 2.62E-01 & 3.37E+06 & 87.1 \\ \hline
        $144^3$ & 7.80E+01 & 3.83E+04 & 8.29E-01 & 3.60E+06 & 94.1 \\ \hline
        $192^3$ & 1.89E+02 & 3.74E+04 & 1.87E+00 & 3.78E+06 & 100.9 \\ \hline
    \end{tabular}
\end{table}

    Tables V and VI show the speedups associated with an MHD code that was run in two
possible configurations. In the first configuration we ran it on the GPU without the {\tt U\_skinny}
trick from Section II, and those results are shown in Table V. In the second configuration we
used the {\tt U\_skinny} trick and ran it on the GPU.

    Table V clearly shows that as we go to higher orders, the speedup suffers. This is because
the {\tt U} variable itself is sent from host to node at the beginning of the timestep and it is sent
back from node to host at the end of the timestep. This turned out to be a very bad choice
because of the performance degradation. We, nevertheless, show it here because it will be very
natural for several CPU-based codes to be written in this fashion and unless the appropriate trick
is used to make a small modification of the code, the GPU advantage will not be realized. Table
VI shows the results from the GPU after the {\tt U\_skinny} trick was implemented correctly. We
see that even on rather small meshes of $96^3$ or $144^3$ zones, a very admirable speedup is achieved.

    We also begin to see another very attractive feature in Table VI, which is that the higher
order schemes show better speedup than the lower order schemes. In other words, a higher order
MHD scheme naturally entails more floating point operations per zone than a lower order
scheme. However, when comparing the third order results to the second order results in Table VI
we see that the number of zones updated per second on the GPU for both schemes is quite
competitive. This is because the third order MHD code is more floating point intensive, with the
result that it utilizes the computational capabilities of the GPU cores much more efficiently. The
similar trend is unfortunately not fully realized when we compare the fourth order results to the
third order results in Table VI; and we will identify the reason in Table VII. The tables also show
that the GPU speedups only become attractive when sufficient work is given to the GPU. For
example, when meshes with $48^3$ zones were given to the GPU, the performance improvement
from the GPU was not quite so dramatic. It is only at $96^3$ or $144^3$ zones per patch that we see
dramatic improvements in the GPU performance relative to the CPU performance in Table VI.
This tells us that rather large patches of MHD mesh need to be processed on a GPU before the
advantage of the GPU is realized.

\begin{table}[!ht]
    \centering
    \caption{comparing the fractional time spent in the ADER subroutine at various
    orders on CPU and GPU for the MHD code.}
    \vspace{2mm}
    \begin{tabular}{|c|c|c|c|c|}
    \hline
              &                 & Fraction of time         &              & Fraction of time         \\
        Zones & Serial time (s) &                  in ADER & GPU time (s) &                  in ADER \\ \hline
        O2 & ~ & ~ & ~ & ~ \\ \hline
        $48^3$ & 6.24E-02 & 0.078 & 5.62E-04 & 0.034 \\ \hline
        $96^3$ & 3.93E-01 & 0.064 & 3.44E-03 & 0.051 \\ \hline
        $144^3$ & 1.22E+00 & 0.058 & 1.06E-02 & 0.05 \\ \hline
        $192^3$ & 2.87E+00 & 0.055 & 2.42E-02 & 0.052 \\ \hline
        ~ & ~ & ~ & ~ & ~ \\ \hline
        O3 & ~ & ~ & ~ & ~ \\ \hline
        $48^3$ & 2.29E-01 & 0.172 & 4.45E-03 & 0.191 \\ \hline
        $96^3$ & 1.50E+00 & 0.143 & 2.76E-02 & 0.272 \\ \hline
        $144^3$ & 4.81E+00 & 0.126 & 8.58E-02 & 0.27 \\ \hline
        $192^3$ & 1.10E+01 & 0.117 & 1.95E-01 & 0.257 \\ \hline
        ~ & ~ & ~ & ~ & ~ \\ \hline
        O4 & ~ & ~ & ~ & ~ \\ \hline
        $48^3$ & 8.11E-01 & 0.262 & 2.77E-02 & 0.537 \\ \hline
        $96^3$ & 4.84E+00 & 0.212 & 1.59E-01 & 0.607 \\ \hline
        $144^3$ & 1.47E+01 & 0.188 & 4.81E-01 & 0.58 \\ \hline
        $192^3$ & 3.37E+01 & 0.178 & 1.07E+00 & 0.573 \\ \hline
    \end{tabular}
\end{table}

    Table VII brings out one of the shortcomings of the GPU. Realize that the caches
available on each GPU core are significantly smaller than the caches available on the CPU cores.
Conceptually, all the MHD codes explored here can be viewed as a WENO reconstruction,
followed by an ADER predictor step and then followed by a Riemann-solver based corrector
step. We were able to write the WENO algorithm in a way that has the smallest possible memory
footprint in the processors cache. We were also able to write the HLL Riemann solver so that it
has the smallest possible memory footprint in the processors cache. These two improvements
caused significant speedups in the CPU and GPU performance of the code. The ADER step that
we used for MHD is based on the ADER algorithm described in Balsara et al. [10], [11] and
already uses the smallest possible arrays that are needed to implement the ADER algorithm
within each zone. The ADER algorithm, therefore, utilizes the CPUs larger cache very
efficiently and is very good at allowing us to make cache-resident calculations on any modern
CPU. However, recall that the GPU cores have significantly smaller caches. As the order of the
ADER scheme increases, the temporary arrays that it needs also utilize more memory. This
memory utilization is small at all orders on a CPU; however it is not small compared to the cache
on a GPU. For that reason, Table VII shows just the ADER performance at various orders on
CPU and GPU. We see that the ADER performance at second order remains excellent on the
GPU, because the second order ADER algorithm can fit inside the GPUs cache. At third order,
Table VII shows a slight degradation in the ADER performance relative to second order. At
fourth order, Table VII shows an even greater degradation in the ADER performance relative to
third order. This is because with increasing order of accuracy the ADER algorithm utilizes the
very tiny GPU cache in an increasingly inefficient way. While this is not very exciting for the
application developer, it highlights a direction in which GPU architectures can be improved in
order for GPUs to become more broadly useful.

    Another way to demonstrate some of the points made in the above paragraph would be to
compare the times for the whole code in Table VI to the time taken for just the ADER update in
Table VII. The fraction of time spent in the ADER predictor step is shown in Table VII. At
second order, and on larger meshes with $96^3$ or $144^3$ zones, we see that the CPU and GPU spend
about the same fraction of time in the ADER predictor step. At third order, and again on larger
meshes with $96^3$ or $144^3$ zones, we see that the fraction of time spend in the ADER step increases
on the GPU compared to the CPU. At fourth order, and again on larger meshes with $96^3$ or $144^3$
zones, we now see that the fraction of time spent in the ADER step is again reasonably constant
for the CPU. However, the fraction of time spent in the ADER step on the GPU increases. It
might be possible to hard-code the ADER algorithm at each order for GPUs, but that would
diminish its portability and malleability of usage. We refer to this as the Predictor bottleneck
because it is currently not possible to arrive at general-purpose predictor algorithms at all orders
where the implementation is small enough to be cache-resident on the very small GPU caches of
present-day GPUs.

\begin{table}[!ht]
    \centering
    \caption{A100 GPU speedup for the MHD code where we exclusively focus on
    the performance differences between the ADER-based scheme and the Runge-Kutta-based scheme.}
    \vspace{2mm}
    \begin{tabular}{|c|c|c|c|c|c|c|}
    \hline
              &           &           &           &           & Relative  &    ~     \\
              &   ADER -  &           &           & Speed up  & Speed up  & Speed up \\
              &    CPU    & RK - CPU  & RK - GPU  & ADER-GPU/ &  RK-GPU/  &  RK-GPU/ \\
        Zones & Zones/sec & Zones/sec & Zones/sec & ADER-CPU  & ADER-CPU  &  RK-CPU  \\ \hline
        O2 & ~ & ~ & ~ & ~ & ~ & ~ \\ \hline
        $48^3$ & 1.38E+05 & 7.53E+04 & 3.42E+06 & 48.2 & 24.77 & 45.42 \\ \hline
        $96^3$ & 1.43E+05 & 7.60E+04 & 7.06E+06 & 92.2 & 49.36 & 92.98 \\ \hline
        $144^3$ & 1.43E+05 & 7.71E+04 & 7.42E+06 & 98.1 & 51.96 & 96.25 \\ \hline
        $192^3$ & 1.34E+05 & 7.33E+04 & 8.12E+06 & 112.9 & 60.42 & 110.8 \\ \hline
        ~ & ~ & ~ & ~ & ~ & ~ & ~ \\ \hline
        O3 & ~ & ~ & ~ & ~ & ~ & ~ \\ \hline
        $48^3$ & 8.29E+04 & 3.28E+04 & 1.97E+06 & 57.2 & 23.7 & 60.0 \\ \hline
        $96^3$ & 8.41E+04 & 3.22E+04 & 3.93E+06 & 103.9 & 46.78 & 122.23 \\ \hline
        $144^3$ & 7.83E+04 & 3.03E+04 & 4.23E+06 & 119.8 & 54.06 & 139.59 \\ \hline
        $192^3$ & 7.50E+04 & 2.78E+04 & 4.39E+06 & 124.3 & 58.55 & 158.12 \\ \hline
        ~ & ~ & ~ & ~ & ~ & ~ & ~ \\ \hline
        O4 & ~ & ~ & ~ & ~ & ~ & ~ \\ \hline
        $48^3$ & 3.58E+04 & 8.96E+03 & 9.29E+05 & 60.0 & 25.97 & 103.69 \\ \hline
        $96^3$ & 3.87E+04 & 9.44E+03 & 1.70E+06 & 87.1 & 43.94 & 180.25 \\ \hline
        $144^3$ & 3.83E+04 & 9.11E+03 & 1.74E+06 & 94.1 & 45.36 & 190.63 \\ \hline
        $192^3$ & 3.74E+04 & 8.88E+03 & 1.77E+06 & 100.9 & 47.2 & 199.12 \\ \hline
    \end{tabular}
\end{table}

    As with Table IV in the previous sub-section, Table VIII in this section considers Runge-
Kutta based schemes. The conclusion from either of these tables is the same. The Runge-Kutta
based MHD code is an inferior performer on GPUs when compared to the ADER based MHD
code. However, the last column of Table VIII mirrors the last column of Table IV. It shows that
as we go to higher orders a Runge-Kutta based code on the GPU will show impressive speedups
telative to the Runge-Kutta based code on the CPU. Furthermore, this speedup improves as we
go to higher order as long as the {\tt U\_skinny} trick is used.

\subsection{CED Results}
    There is an essential, and structural, difference between CFD and MHD codes on the one
hand and CED codes on the other hand. For CFD and MHD codes, the WENO based
reconstruction step, the ADER based predictor step and the Riemann solver based corrector step
all require roughly comparable amounts of floating point operations. In a CED code, the WENO
based reconstruction step has to be entirely divergence-preserving (Balsara [9] or Balsara et al.
[16]) which makes it very inexpensive. Furthermore, the multidimensional Riemann solver used
in the corrector step of a CED code is extremely inexpensive (see eqns. (5.5) and (5.6) from
Balsara et al. [15]). As a result, the reconstruction and the corrector steps in a CED code are
unusually lightweight and require very few floating point operations. However, realize too that
most CED applications do have unusually large conductivities, which results in very stiff source
terms. ADER schemes that treat stiff source terms (Dumbser et al. [23], Balsara et al. [16]) have
been designed, but they are considerably more memory and floating point intensive than ADER
schemes for treating problems without stiff source terms. While diagonally implicit Runge-Kutta
schemes have also been used for treating stiff source terms in CED (Balsara et al. [13]), they
tend to be even more inefficient than ADER schemes when used in finite volume context for
CED applications. (This is why we dont show any Runge-Kutta based results for CED.) Thus
CED is a harsher algorithmic terrain for GPU computing compared to CFD and MHD.

\begin{table}[!ht]
    \caption{A100 GPU speedup for the CED code when the {\tt U\_skinny} trick was not used.}
    \vspace{2mm}
    \centering
    \begin{tabular}{|c|c|c|c|c|c|c|c|}
    \hline
              &                 & Zones /         &              & Zones /      &          \\ 
        Zones & Serial time (s) & Serial time (s) & GPU time (s) & GPU time (s) & Speed up \\ \hline
        O2 & ~ & ~ & ~ & ~ & ~ \\ \hline
        $48^3$ & 4.09E-01 & 2.71E+05 & 1.20E-02 & 9.23E+06 & 34.1 \\ \hline
        $96^3$ & 2.69E+00 & 3.29E+05 & 7.63E-02 & 1.16E+07 & 35.2 \\ \hline
        $144^3$ & 9.94E+00 & 3.00E+05 & 2.37E-01 & 1.26E+07 & 41.9 \\ \hline
        $192^3$ & 2.17E+01 & 3.26E+05 & 5.21E-01 & 1.36E+07 & 41.8 \\ \hline
        ~ & ~ & ~ & ~ & ~ & ~ \\ \hline
        O3 & ~ & ~ & ~ & ~ & ~ \\ \hline
        $48^3$ & 1.19E+00 & 9.31E+04 & 3.41E-02 & 3.24E+06 & 34.8 \\ \hline
        $96^3$ & 8.19E+00 & 1.08E+05 & 2.14E-01 & 4.14E+06 & 38.3 \\ \hline
        $144^3$ & 2.64E+01 & 1.13E+05 & 6.73E-01 & 4.44E+06 & 39.2 \\ \hline
        $192^3$ & 6.41E+01 & 1.10E+05 & 1.46E+00 & 4.84E+06 & 43.8 \\ \hline
        ~ & ~ & ~ & ~ & ~ & ~ \\ \hline
        O4 & ~ & ~ & ~ & ~ & ~ \\ \hline
        $48^3$ & 5.74E+00 & 1.93E+04 & 1.24E-01 & 8.88E+05 & 46.1 \\ \hline
        $96^3$ & 3.72E+01 & 2.38E+04 & 7.15E-01 & 1.24E+06 & 52 \\ \hline
        $144^3$ & 1.12E+02 & 2.67E+04 & 2.13E+00 & 1.41E+06 & 52.7 \\ \hline
        $180^3$ & 2.08E+02 & 2.81E+04 & 3.98E+00 & 1.46E+06 & 52.2 \\ \hline
    \end{tabular}
\end{table}

\begin{table}[!ht]
    \caption{A100 GPU speedup for the CED code when the {\tt U\_skinny} trick was used.}
    \vspace{2mm}
    \centering
    \begin{tabular}{|c|c|c|c|c|c|c|c|}
    \hline
              &                 & Zones /         &              & Zones /      &          \\ 
        Zones & Serial time (s) & Serial time (s) & GPU time (s) & GPU time (s) & Speed up \\ \hline
        O2 & ~ & ~ & ~ & ~ & ~ \\ \hline
        $48^3$ & 6.28E-01 & 1.76E+05 & 1.29E-02 & 8.55E+06 & 48.54 \\ \hline
        $96^3$ & 4.56E+00 & 1.94E+05 & 7.60E-02 & 1.16E+07 & 60.02 \\ \hline
        $144^3$ & 1.40E+01 & 2.14E+05 & 2.50E-01 & 1.20E+07 & 55.95 \\ \hline
        $192^3$ & 3.18E+01 & 2.22E+05 & 5.69E-01 & 1.24E+07 & 55.93 \\ \hline
        ~ & ~ & ~ & ~ & ~ & ~ \\ \hline
        O3 & ~ & ~ & ~ & ~ & ~ \\ \hline
        $48^3$ & 1.90E+00 & 5.83E+04 & 3.27E-02 & 3.38E+06 & 58.01 \\ \hline
        $96^3$ & 1.29E+01 & 6.88E+04 & 2.01E-01 & 4.41E+06 & 64.04 \\ \hline
        $144^3$ & 4.20E+01 & 7.11E+04 & 6.38E-01 & 4.68E+06 & 65.85 \\ \hline
        $192^3$ & 9.20E+01 & 7.69E+04 & 1.46E+00 & 4.85E+06 & 63.12 \\ \hline
        ~ & ~ & ~ & ~ & ~ & ~ \\ \hline
        O4 & ~ & ~ & ~ & ~ & ~ \\ \hline
        $48^3$ & 7.28E+00 & 1.52E+04 & 9.21E-02 & 1.20E+06 & 78.98 \\ \hline
        $96^3$ & 4.72E+01 & 1.88E+04 & 5.72E-01 & 1.55E+06 & 82.44 \\ \hline
        $144^3$ & 1.42E+02 & 2.10E+04 & 1.79E+00 & 1.66E+06 & 79.2 \\ \hline
        $180^3$ & 2.67E+02 & 2.18E+04 & 3.42E+00 & 1.71E+06 & 78.23 \\ \hline
    \end{tabular}
\end{table}
    Tables IX and X show the performance of CED codes when the {\tt U\_skinny} trick is not
used and when it is used. We see from Table X that the trick does have a rather positive effect on
the overall speedup. However, the speedups are not as impressive as the ones that we have seen
in prior sections for CFD and MHD. This is understandable considering that the ADER predictor
step is so dominant for CED codes, as we will see in the next paragraph.

\begin{table}[!ht]
    \centering
    \caption{comparing the fractional time spent in the ADER subroutine at various
    orders on CPU and GPU for the CED code.}
    \vspace{2mm}
    \begin{tabular}{|c|c|c|c|c|}
    \hline
              &                 & Fraction of time         &              & Fraction of time         \\
        Zones & Serial time (s) &                  in ADER & GPU time (s) &                  in ADER \\ \hline
        O2 & ~ & ~ & ~ & ~ \\ \hline
        $48^3$ & 3.18E-01 & 0.506 & 6.95E-03 & 0.538 \\ \hline
        $96^3$ & 2.09E+00 & 0.458 & 3.56E-02 & 0.468 \\ \hline
        $144^3$ & 6.63E+00 & 0.475 & 1.09E-01 & 0.437 \\ \hline
        $192^3$ & 1.45E+01 & 0.457 & 2.46E-01 & 0.433 \\ \hline
        ~ & ~ & ~ & ~ & ~ \\ \hline
        O3 & ~ & ~ & ~ & ~ \\ \hline
        $48^3$ & 1.48E+00 & 0.778 & 2.41E-02 & 0.738 \\ \hline
        $96^3$ & 9.59E+00 & 0.746 & 1.49E-01 & 0.742 \\ \hline
        $144^3$ & 2.95E+01 & 0.702 & 4.61E-01 & 0.723 \\ \hline
        $192^3$ & 6.39E+01 & 0.694 & 1.06E+00 & 0.724 \\ \hline
        ~ & ~ & ~ & ~ & ~ \\ \hline
        O4 & ~ & ~ & ~ & ~ \\ \hline
        $48^3$ & 6.34E+00 & 0.871 & 7.58E-02 & 0.823 \\ \hline
        $96^3$ & 3.98E+01 & 0.844 & 4.71E-01 & 0.823 \\ \hline
        $144^3$ & 1.18E+02 & 0.828 & 1.47E+00 & 0.817 \\ \hline
        $180^3$ & 2.21E+02 & 0.826 & 2.78E+00 & 0.813 \\ \hline
    \end{tabular}
\end{table}

    Table XI focuses on the fractional time taken by the ADER predictor step for CED on
CPUs and GPUs at various orders. It gives us greater clarity of perspective. At second order, the
ADER predictor step is still quite light weight and takes up comparable fractions of time on the
CPU and the GPU. In other words, at second order, it is possible to arrive at cache-friendly
versions of the ADER predictor step on both the CPU and the GPU. However, at third order, we
see that the fraction of time taken in the ADER predictor step is considerably larger on the GPU
than on the CPU. This tells us that at third order and on CPUs, the ADER can still function in a
relatively cache-friendly mode. However, at third order and on GPUs, the ADER is not as
efficient because it is not so cache-friendly. The fourth order results in Table XI are even more
interesting. We see that the ADER seems to entirely dominate the time taken for a timestep. This
is so on CPUs and GPUs, which tells us that in either computational environment the ADER is
not so cache-friendly. Even though this result is not as optimistic as one would desire, it reveals
an opportunity. It shows us that there might be room for improving the ADER algorithm at
higher orders by finding more cache-friendly variants. It also tells us that machines with larger
caches would be more suitable for this type of application.

\section{Accuracy Analysis CPU v/s GPU}
The astute reader might still point out that the new thing that makes GPUs competitive with
CPUs is indeed the inclusion of error correction in the GPUs. Therefore, we are interested in
providing some perspective on how good the error correction indeed is on the GPUs. This is
done by taking standard test problems that are used for accuracy analysis and running them on a
CPU and gathering the results. We then run the same test for exactly the same number of
timesteps and identical stopping time on a GPU. Sub-section IV.1 documents CFD results; Sub-
section IV.2 documents MHD results and Sub-section IV.3 documents results from CED. The
happy conclusion from these sections is that the accuracy numbers are virtually identical up to
about nine or ten digits of accuracy. As a result, the reader will find the tables that show accuracy
below virtually identical for CPU and GPU.

\subsection{CFD Results}
    The result shown is from a three-dimensional version of the isentropic hydrodynamic
vortex test problem from Balsara and Shu [9]; see Sub-section IV.b of that work. Table XII
shows the $L_1$ and $L_\infty$ errors at second, third and fourth orders for our CFD test problem. We also
show that the schemes on both CPU and GPU do indeed reach their target accuracies. Up to four
significant digits, the error from the CPU and GPU look identical. As a result, the last column
gives the absolute value of the difference between the error as measured on the CPU and the
GPU. This drives home the point that the difference between the CPU result and the result from
the GPU is indeed in the digits that would be considered insignificant for scientific computation.
The table shows that on a CPU, as well as on a GPU, virtually identical accuracies are obtained.

\begin{table}[!ht]
    \centering
    \caption{Errors at second, third and fourth orders for our CFD 
    test problem. We show the error in the density variable for the hydrodynamical 
    vortex. Up to four significant digits, the error from the CPU and GPU look 
    identical. As a result, the last column gives the absolute value of the 
    difference between the error as measured on the CPU and the GPU.}
    \vspace{2mm}
    \begin{tabular}{|c|c|c|c|c|}
    \hline
 CPU v/s GPU &                 &                 &                & $|$ CPU $L_1$ Error - \\ 
 O2; Zones   & $L_1$ Error CPU & $L_1$ Error GPU & $L_1$ Accuracy &     GPU $L_1$ Error $|$  \\ \hline
        $48^3$ & 1.60E-03 & 1.60E-03 & -- & 1.98E-15 \\ \hline
        $96^3$ & 4.03E-04 & 4.03E-04 & 1.99 & 6.10E-17 \\ \hline
        $144^3$ & 1.63E-04 & 1.63E-04 & 2.23 & 2.20E-17 \\ \hline
        $192^3$ & 8.74E-05 & 8.74E-05 & 2.16 & 1.31E-10 \\ \hline
CPU v/s GPU &                 &                 &                              & $|$ CPU $L_\infty$ Error - \\ 
O2; Zones   & $L_\infty$ Error CPU & $L_\infty$ Error GPU & $L_\infty$ Accuracy &    GPU $L_\infty$ Error $|$\\\hline
        $48^3$ & 2.82E-02 & 2.82E-02 & -- & 1.23E-14 \\ \hline
        $96^3$ & 7.47E-03 & 7.47E-03 & 1.92 & 1.10E-16 \\ \hline
        $144^3$ & 3.00E-03 & 3.00E-03 & 2.25 & 4.50E-16 \\ \hline
        $192^3$ & 1.55E-03 & 1.55E-03 & 2.29 & 3.10E-09 \\ \hline
       CPU v/s GPU &                 &                 &                & $|$ CPU $L_1$ Error - \\ 
       O3; Zones   & $L_1$ Error CPU & $L_1$ Error GPU & $L_1$ Accuracy &  GPU $L_1$ Error $|$  \\ \hline
        $48^3$ & 1.09E-03 & 1.09E-03 & -- & 3.90E-16 \\ \hline
        $96^3$ & 2.25E-04 & 2.25E-04 & 2.27 & 9.76E-19 \\ \hline
        $144^3$ & 7.07E-05 & 7.07E-05 & 2.85 & 8.78E-17 \\ \hline
        $192^3$ & 3.03E-05 & 3.03E-05 & 2.94 & 3.20E-17 \\ \hline
      CPU v/s GPU &                 &                 &                & $|$ CPU $L_\infty$ Error - \\ 
      O3; Zones   & $L_\infty$ Error CPU & $L_\infty$ Error GPU & $L_\infty$ Accuracy & GPU $L_\infty$ Error $|$ \\\hline
        $48^3$ & 1.84E-02 & 1.84E-02 & -- & 6.00E-16 \\ \hline
        $96^3$ & 3.21E-03 & 3.21E-03 & 2.52 & 5.60E-16 \\ \hline
        $144^3$ & 1.06E-03 & 1.06E-03 & 2.74 & 1.66E-15 \\ \hline
        $192^3$ & 4.63E-04 & 4.63E-04 & 2.87 & 3.33E-16 \\ \hline
       CPU v/s GPU &                 &                 &                & $|$ CPU $L_1$ Error - \\ 
       O4; Zones   & $L_1$ Error CPU & $L_1$ Error GPU & $L_1$ Accuracy &  GPU $L_1$ Error $|$  \\ \hline
        $48^3$ & 1.03E-03 & 1.03E-03 & -- & 1.20E-14 \\ \hline
        $96^3$ & 9.46E-05 & 9.46E-05 & 3.44 & 8.48E-16 \\ \hline
        $144^3$ & 1.84E-05 & 1.84E-05 & 4.04 & 2.07E-17 \\ \hline
        $192^3$ & 5.48E-06 & 5.48E-06 & 4.22 & 2.16E-17 \\ \hline
      CPU v/s GPU &                 &                 &                & $|$ CPU $L_\infty$ Error - \\ 
      O4; Zones   & $L_\infty$ Error CPU & $L_\infty$ Error GPU & $L_\infty$ Accuracy & GPU $L_\infty$ Error $|$\\\hline
        $48^3$ & 6.50E-02 & 6.50E-02 & -- & 5.35E-13 \\ \hline
        $96^3$ & 1.15E-02 & 1.15E-02 & 2.5 & 1.21E-14 \\ \hline
        $144^3$ & 9.64E-04 & 9.64E-04 & 6.11 & 1.33E-15 \\ \hline
        $192^3$ & 4.38E-04 & 4.38E-04 & 2.74 & 4.44E-16 \\ \hline
    \end{tabular}
\end{table}

Table XII shows the $L_1$ and $L_\infty$ errors at second, third and fourth orders for our CFD test
problem. We show the error in the density variable for the hydrodynamical vortex. Up to
four significant digits, the error from the CPU and GPU look identical. As a result, the last
column gives the absolute value of the difference between the error as measured on the
CPU and the GPU.

\subsection{MHD Results}
    The result shown is from a three-dimensional version of the MHD vortex test problem
from Balsara [3]; see Section 6 of that work. Table XIII shows the $L_1$ and $L_\infty$ errors at second,
third and fourth orders for our MHD test problem. The results are presented in the same style as
in Table XII. The table shows that on a CPU, as well as on a GPU, virtually identical accuracies
are obtained.

\begin{table}[!ht]
    \centering
    \caption{Errors at second, third and fourth orders for our MHD test problem. 
    We show the error in the x-component of the magnetic variable for the MHD vortex.
    The results are presented in the same style as in Table XII.}
    \vspace{2mm}
    \begin{tabular}{|c|c|c|c|c|}
    \hline
 CPU v/s GPU &                 &                 &                & $|$ CPU $L_1$ Error - \\ 
 O2; Zones   & $L_1$ Error CPU & $L_1$ Error GPU & $L_1$ Accuracy &     GPU $L_1$ Error $|$  \\ \hline
        $48^3$ & 4.65E-03 & 4.65E-03 & -- & 1.04E-17 \\ \hline
        $96^3$ & 1.24E-03 & 1.24E-03 & 1.91 & 9.97E-18 \\ \hline
        $144^3$ & 5.58E-04 & 5.58E-04 & 1.97 & 4.99E-18 \\ \hline
        $192^3$ & 3.15E-04 & 3.15E-04 & 1.99 & 6.60E-17 \\ \hline
CPU v/s GPU &                 &                 &                              & $|$ CPU $L_\infty$ Error - \\ 
O2; Zones   & $L_\infty$ Error CPU & $L_\infty$ Error GPU & $L_\infty$ Accuracy &    GPU $L_\infty$ Error $|$\\\hline
        $48^3$ & 6.18E-02 & 6.18E-02 & -- & 1.66E-14 \\ \hline
        $96^3$ & 1.54E-02 & 1.54E-02 & 2.01 & 3.30E-15 \\ \hline
        $144^3$ & 6.83E-03 & 6.83E-03 & 2 & 1.19E-15 \\ \hline
        $192^3$ & 3.85E-03 & 3.85E-03 & 1.99 & 4.22E-15 \\ \hline
 CPU v/s GPU &                 &                 &                & $|$ CPU $L_1$ Error - \\ 
 O3; Zones   & $L_1$ Error CPU & $L_1$ Error GPU & $L_1$ Accuracy &     GPU $L_1$ Error $|$  \\ \hline
        $48^3$ & 1.43E-03 & 1.43E-03 & -- & 1.18E-15 \\ \hline
        $96^3$ & 1.92E-04 & 1.92E-04 & 2.9 & 8.20E-17 \\ \hline
        $144^3$ & 5.73E-05 & 5.73E-05 & 2.98 & 1.68E-17 \\ \hline
        $192^3$ & 2.42E-05 & 2.42E-05 & 2.99 & 3.99E-17 \\ \hline
CPU v/s GPU &                 &                 &                              & $|$ CPU $L_\infty$ Error - \\ 
O3; Zones   & $L_\infty$ Error CPU & $L_\infty$ Error GPU & $L_\infty$ Accuracy &    GPU $L_\infty$ Error $|$\\\hline
        $48^3$ & 2.46E-02 & 2.46E-02 & -- & 1.50E-15 \\ \hline
        $96^3$ & 3.01E-03 & 3.01E-03 & 3.03 & 3.90E-16 \\ \hline
        $144^3$ & 8.32E-04 & 8.32E-04 & 3.17 & 2.16E-15 \\ \hline
        $192^3$ & 3.59E-04 & 3.59E-04 & 2.92 & 8.33E-16 \\ \hline
 CPU v/s GPU &                 &                 &                & $|$ CPU $L_1$ Error - \\ 
 O4; Zones   & $L_1$ Error CPU & $L_1$ Error GPU & $L_1$ Accuracy &     GPU $L_1$ Error $|$  \\ \hline
        $48^3$ & 1.40E-03 & 1.40E-03 & -- & 7.56E-14 \\ \hline
        $96^3$ & 6.16E-05 & 6.16E-05 & 4.51 & 5.69E-18 \\ \hline
        $144^3$ & 1.02E-05 & 1.02E-05 & 4.43 & 2.31E-17 \\ \hline
        $192^3$ & 3.00E-06 & 3.00E-06 & 4.26 & 1.98E-17 \\ \hline
CPU v/s GPU &                 &                 &                              & $|$ CPU $L_\infty$ Error - \\ 
O4; Zones   & $L_\infty$ Error CPU & $L_\infty$ Error GPU & $L_\infty$ Accuracy &    GPU $L_\infty$ Error $|$\\\hline
        $48^3$ & 9.32E-02 & 9.32E-02 & -- & 1.03E-12 \\ \hline
        $96^3$ & 4.45E-03 & 4.45E-03 & 4.39 & 3.33E-15 \\ \hline
        $144^3$ & 6.57E-04 & 6.57E-04 & 4.72 & 1.67E-15 \\ \hline
        $192^3$ & 1.72E-04 & 1.72E-04 & 4.66 & 6.33E-15 \\ \hline
    \end{tabular}
\end{table}

Table XIII shows the $L_1$ and $L_\infty$ errors at second, third and fourth orders for our MHD test
problem. We show the error in the x-component of the magnetic variable for the MHD
vortex. The results are presented in the same style as in Table XII.

\subsection{CED Results}
    The result shown is from a three-dimensional version of the electromagnetic wave
propagation test problem from Balsara et al. [16]; see Sub-section 5.1 of that work. Table XIV
shows the $L_1$ and $L_\infty$ errors at second, third and fourth orders for our CED test problem. The
results are presented in the same style as in Table XII. The table shows that on a CPU, as well as
on a GPU, virtually identical accuracies are obtained.

\begin{table}[!ht]
    \centering
    \caption{Errors at second, third and fourth orders for our CED test problem. 
    We show the error in the y-component of the electric displacement variable 
    for the electromagnetic wave. The results are presented in the same style as in Table XII.}
    \vspace{2mm}
    \begin{tabular}{|c|c|c|c|c|}
    \hline
 CPU v/s GPU &                 &                 &                & $|$ CPU $L_1$ Error - \\ 
 O2; Zones   & $L_1$ Error CPU & $L_1$ Error GPU & $L_1$ Accuracy &     GPU $L_1$ Error $|$  \\ \hline
        $16^3$ & 1.03E-04 & 1.03E-04 & -- & 1.99E-18 \\ \hline
        $32^3$ & 2.35E-05 & 2.35E-05 & 2.13 & 0.00E+00 \\ \hline
        $64^3$ & 5.68E-06 & 5.68E-06 & 2.05 & 2.64E-18 \\ \hline
        $128^3$ & 1.41E-06 & 1.41E-06 & 2.01 & 1.90E-18 \\ \hline
CPU v/s GPU &                 &                 &                              & $|$ CPU $L_\infty$ Error - \\ 
O2; Zones   & $L_\infty$ Error CPU & $L_\infty$ Error GPU & $L_\infty$ Accuracy &    GPU $L_\infty$ Error $|$\\\hline
        $16^3$ & 1.63E-04 & 1.63E-04 & -- & 0.00E+00 \\ \hline
        $32^3$ & 3.67E-05 & 3.67E-05 & 2.15 & 2.00E-18 \\ \hline
        $64^3$ & 8.92E-06 & 8.92E-06 & 2.04 & 2.65E-18 \\ \hline
        $128^3$ & 2.21E-06 & 2.21E-06 & 2.01 & 3.14E-18 \\ \hline
 CPU v/s GPU &                 &                 &                & $|$ CPU $L_1$ Error - \\ 
 O3; Zones   & $L_1$ Error CPU & $L_1$ Error GPU & $L_1$ Accuracy &     GPU $L_1$ Error $|$  \\ \hline
        $16^3$ & 3.59E-05 & 3.59E-05 & -- & 0.00E+00 \\ \hline
        $32^3$ & 4.39E-06 & 4.39E-06 & 3.03 & 3.05E-20 \\ \hline
        $64^3$ & 5.41E-07 & 5.41E-07 & 3.02 & 1.11E-19 \\ \hline
        $128^3$ & 6.71E-08 & 6.71E-08 & 3.01 & 2.33E-19 \\ \hline
CPU v/s GPU &                 &                 &                              & $|$ CPU $L_\infty$ Error - \\ 
O3; Zones   & $L_\infty$ Error CPU & $L_\infty$ Error GPU & $L_\infty$ Accuracy &    GPU $L_\infty$ Error $|$\\\hline
        $16^3$ & 5.47E-05 & 5.47E-05 & -- & 1.50E-18 \\ \hline
        $32^3$ & 6.78E-06 & 6.78E-06 & 3.01 & 0.00E+00 \\ \hline
        $64^3$ & 8.43E-07 & 8.43E-07 & 3.01 & 1.09E-18 \\ \hline
        $128^3$ & 1.05E-07 & 1.05E-07 & 3.01 & 4.34E-18 \\ \hline
 CPU v/s GPU &                 &                 &                & $|$ CPU $L_1$ Error - \\ 
 O4; Zones   & $L_1$ Error CPU & $L_1$ Error GPU & $L_1$ Accuracy &     GPU $L_1$ Error $|$  \\ \hline
        $16^3$ & 2.04E-06 & 2.04E-06 & -- & 3.90E-19 \\ \hline
        $32^3$ & 1.14E-07 & 1.14E-07 & 4.16 & 6.90E-20 \\ \hline
        $64^3$ & 7.09E-09 & 7.09E-09 & 4.01 & 2.43E-19 \\ \hline
        $128^3$ & 4.52E-10 & 4.52E-10 & 3.97 & 6.04E-19 \\ \hline
CPU v/s GPU &                 &                 &                              & $|$ CPU $L_\infty$ Error - \\ 
O4; Zones   & $L_\infty$ Error CPU & $L_\infty$ Error GPU & $L_\infty$ Accuracy &    GPU $L_\infty$ Error $|$\\\hline
        $16^3$ & 3.39E-06 & 3.39E-06 & -- & 2.60E-18 \\ \hline
        $32^3$ & 1.71E-07 & 1.71E-07 & 4.31 & 2.28E-18 \\ \hline
        $64^3$ & 1.12E-08 & 1.12E-08 & 3.93 & 2.17E-18 \\ \hline
        $128^3$ & 7.05E-10 & 7.05E-10 & 3.99 & 1.36E-20 \\ \hline
    \end{tabular}
\end{table}

Table XIV shows the $L_1$ and $L_\infty$ errors at second, third and fourth orders for our CED test
problem. We show the error in the y-component of the electric displacement variable for
the electromagnetic wave. The results are presented in the same style as in Table XII.

\section{Scalability Study with Comparable Numbers of GPUs and CPUs}
\begin{figure}[htp]
    \includegraphics[scale=0.25]{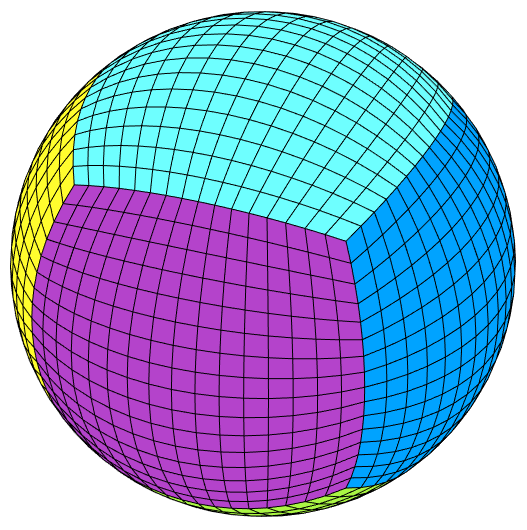} \hspace{4mm}
    \includegraphics[scale=0.17]{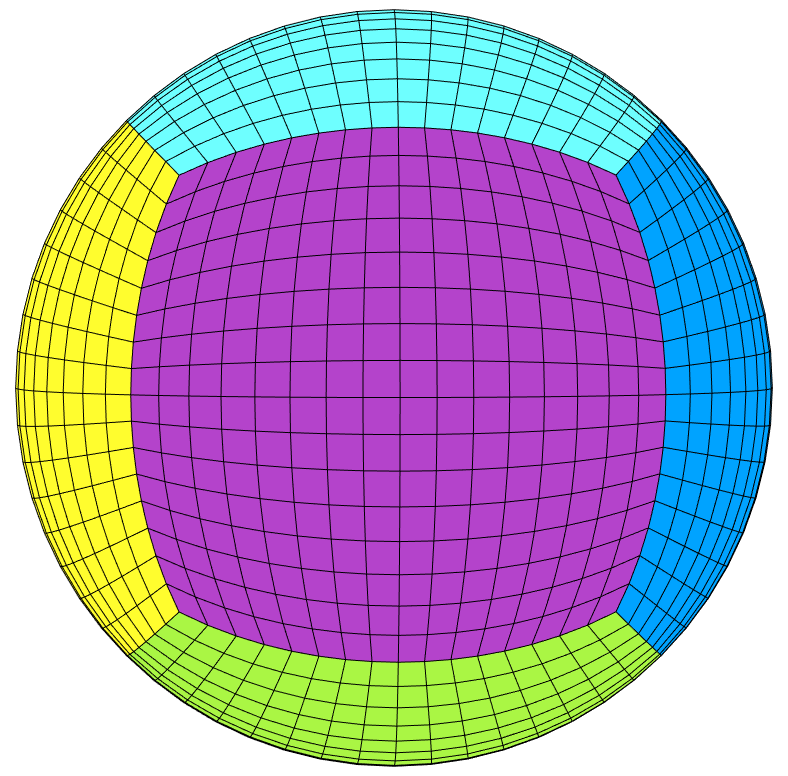} \hspace{4mm}
    \includegraphics[scale=0.25]{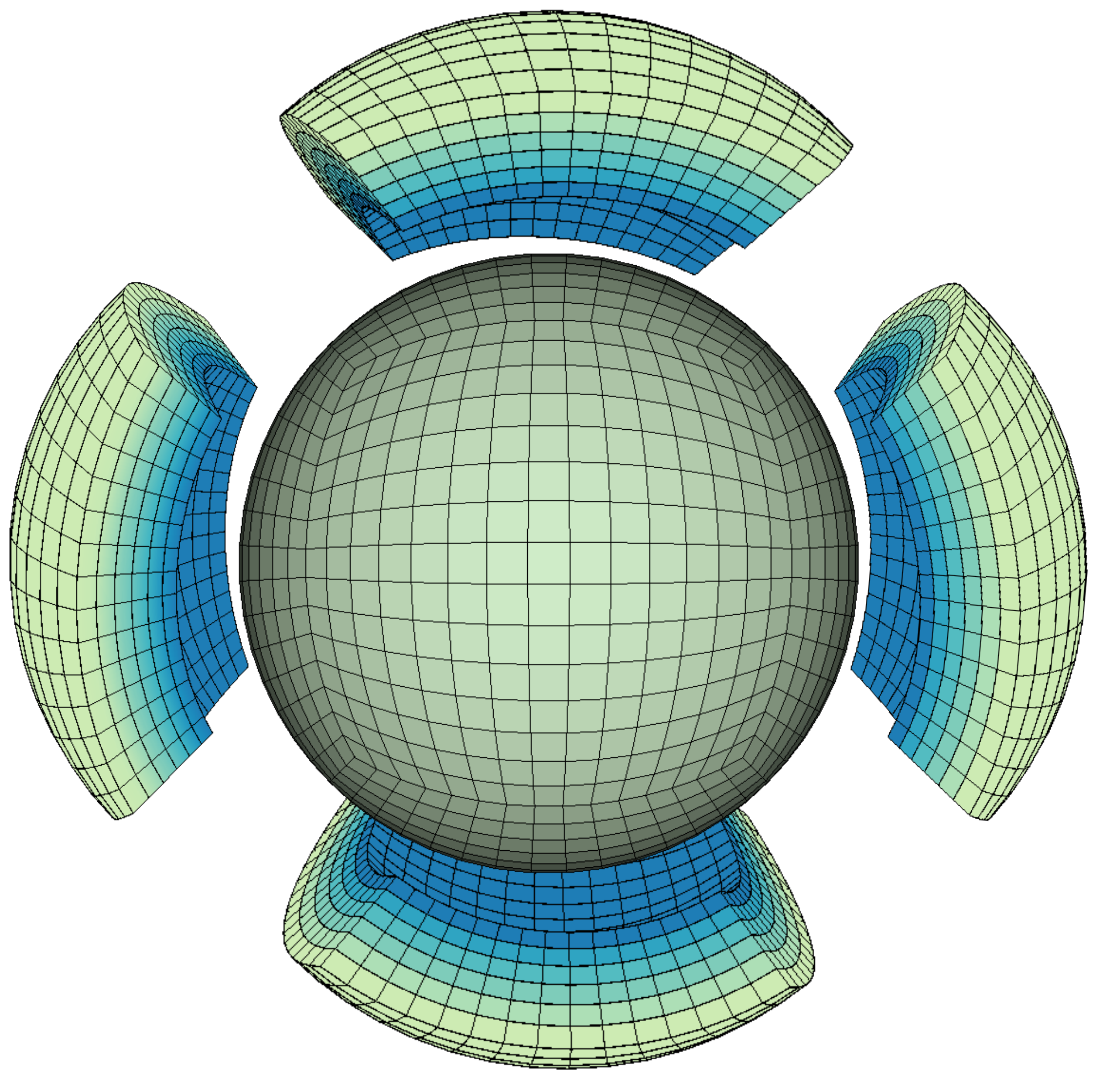}\par
    \hspace{20mm}(a)\hspace{48mm}(b)\hspace{52mm}(c)\par
\caption{ Meshing of an actual physical problem.
    (a) A 2D cube sphere mesh is constructed on the surface of a unit sphere
    (b) Rotated version of the same 2D cube sphere mesh.
    (c) The 2D surface mesh in (a) and (b) results in six sectors, and five of
    those sectors which face the reader are shown with different colors.
    The problem consists of doing 3D MHD calculations on the exterior volume of
    the sphere. To that end, the 2D surface mesh is extended in the third, i.e.
    radial, direction. Also (c) only shows the 3D mesh associated with
    four of the possible six sectors. The resulting mesh gives us a uniform
    coverage of the 3D volume that is exterior to the sphere and is known as a
    3D cube sphere mesh }
\label{fig:cs_mesh}
\end{figure}

    Section II has shown us the different types of techniques, tricks and algorithms needed
for extracting performance from a GPU; and Section III has shown us that the ideas that we have
developed indeed work in idealized situations. We now show how well these ideas work in an
end to end scientific application. It is very illustrative to show a real-world application problem
and carry out a scalability study involving GPUs and comparable number of multicore CPUs for
such a problem. The problem we use is from astrophysical applications, see Subramanian et al.
[37]. The only difference between the cited research and the current application is that now we
wish to do the problem with more physics, which calls for a larger mesh and a more
sophisticated mesh configuration. Fig. 8 shows us the meshing of an actual problem. Of course
the number of zones have been down-sampled in the figure to make it easier for the reader to
visualize the geometry of the problem. Here the physical problem that we want to solve consists
of simulating fully 3D MHD in the volume that is exterior to a sphere, as shown in Fig. 8. Figs.
8a and 8b show us different views of the 2D surface of a unit sphere, showing us that the surface
can be meshed with six sectors. The five sectors that face the reader are shown with different
colors; the sixth sector is hidden from view. Fig. 8c shows us how the 2D sectors can be
converted into 3D sectors by extending the angular mesh from the spherical surface in the radial
direction (i.e. the third direction). Only four of the possible six 3D sectors are shown in Fig. 8c.
In Fig. 8 we show only a few zones in each direction, but an actual computation involves
hundreds of zones in each direction of a 3D sector. Such a 3D mesh that covers the volume that
is exterior to a sphere is known as a cube sphere mesh. The goal is to carry out a 3D MHD
simulation on such a cube sphere mesh. The scalability study will be dictated by the real-world
physics application that is to be done on a GPU-rich machine.

    The physics of the problem is as follows. We need to cover the sphere with an angular
resolution of less than one degree. To that end, we decided that each sector should have
$100 \times 100$ zones in each of the two angular directions on the unit sphere. Furthermore, in the
radial direction the physics of the problem is such that we need 400 zones in the radial direction.
Therefore, each sector of the type shown in Fig. 8c should have $100 \times 100 \times 400$ zones; and please
realize that we have six such sectors in Fig. 8c. On evaluating the application requirements, it
became evident to us that all the data associated with $100 \times 100 \times 400$ zones would not fit within
the RAM memory of a single GPU. Therefore, the data would have to be chunked into smaller
patches and those patches would have to be sent to the device for update and the minimal
updated data would have to be brought back to the host. We found that $100 \times 100 \times 100$ zone
patches could fit on each of the GPUs that we had access to. As a result, we would have to start
our scalability study with 4 such patches per sector and 6 sectors, resulting in 24 patches (each
patch having $100 \times 100 \times 100$ zones) for the full spherical calculation. We also found that the host
could accommodate 8 patches on each of the nodes of the supercomputer we were using. This
paragraph serves to show us that the hardware limitations, and the physics of the problem,
strongly limit how a calculation can be done on a GPU-rich supercomputer. We discuss the
details of the machine that we were using next.

    We used the Delta machine at NCSA along with the NVIDIA NVHPC (version 22.2)
Fortran compiler. Delta has various configurations. The configuration that we had access to had
100 nodes. Each node had one AMD Milan 64 core CPU (the host) and four A100 GPUs with 40
GB HBM2 RAM and NVLink (the devices). We wished to compare the performance of a certain
number of GPUs with the same number of multicore CPUs. As a result, given the size of the
physical problem, it was felt that we could do a weak scalability study that went from 12 GPUs
to 96 GPUs and the same could be compared with a comparable scalability study that went from
12 CPUs to 96 CPUs. Each node of Delta had sufficient main RAM to hold 8 patches with
$100 \times 100 \times 100$ zones per patch. Delta also has only 100 nodes, each with one CPU, so the
scalability study would be restricted to a maximum of 96 CPUs, and therefore 96 GPUs.

    The application demands high angular resolution. As a result, for the 12 GPU run, we
used $100^3$ zones per patch and 24 patches. On the CPUs, we have two options, as follows:-
\begin{itemize}
    \item This first option consists of using only one-sided MPI-3 messaging. To carry out the identical
run on 12 CPUs (with 64 cores per CPU), we had to use $25 \times 25 \times 50$ zones per patch and 768
patches. Please check that the number of zones is the same; i.e., $24 \times 100^3$ zones on the 12 GPUs
is same as $25 \times 25 \times 50 \times 12 \times 64$ zones on the 12 CPUs with 64 cores per CPU. With each
doubling of CPUs or GPUs, we just doubled the number of patches. One-sided MPI-3
communication was used for the messaging across patches in order to have a messaging strategy
with the lowest latency and highest bandwidth. (Please see Garain, Balsara and Reid [26] to
realize that the MPI-3 messaging strategy is optimal and that the on-processor performance of
the code is also optimal for this category of code.) All the MPI communication had to be done
through the CPUs because the GPUs never had enough RAM to concurrently hold all the data for
the whole problem.
    \item The second option consists of using a hybrid parallelization model on the CPUs based on
using 64-way OpenMP-based parallelization within a CPU and MPI-3 parallelization across
CPUs. This strategy allows us to maintain identical patch sizes on GPUs and CPUs. The CPU
patches now had $100 \times 100 \times 100$ zones; i.e., the same as on the GPUs. The benefit of this strategy
is that it allows much fewer MPI messages to be exchanged, with each message being larger. The
downside of this strategy is that one has to include OpenMP parallelization in addition to MPI
parallelization. In this approach too, all the MPI communication had to be done through the
CPUs because the GPUs never had enough RAM to concurrently hold all the data for the whole
problem.
\end{itemize}

The discussion in this paragraph is intended to bring the reader face to face with the harsh
realities and substantial limitations that a computationalist faces when simulating a large real-
word scientific application on an actual GPU-equipped supercomputer!

\begin{table}[!ht]
    \caption{Scalability data when all the computations were done on CPUs with MPI-3 messaging 
    between CPUs}
    \vspace{2mm}
    \centering
    \begin{tabular}{ c c c c c }
    \hline
         CPU & & & & \\ 
   Parallelization & 64 Core  Processors & Zones & CPU time & CPU  Zones/sec \\ \hline
                             ~ & 12 & 24000000  & 2.45E+00 & 9.79E+06 \\   
           MPI Only          ~ & 24 & 48000000  & 3.78E+00 & 1.27E+07 \\   
                             ~ & 48 & 96000000  & 6.46E+00 & 1.49E+07 \\   
                             ~ & 96 & 192000000 & 1.14E+01 & 1.68E+07 \\    \hline
    \end{tabular}
\end{table}

\begin{table}[!ht]
    \caption{Scalability data when all the computations were done on CPUs with MPI-3 messaging
    between CPUs}
    \vspace{2mm}
    \centering
    \begin{tabular}{ c c c c c }
    \hline
         CPU & & & & \\ 
  Parallelization & 64 Core  Processors & Zones & CPU time & CPU  Zones/sec \\ \hline
                               ~ & 12 & 24000000  & 3.25E+00 & 7.39E+06 \\  
          OpenMP + MPI         ~ & 24 & 48000000  & 3.62E+00 & 1.33E+07 \\  
            Hybrid             ~ & 48 & 96000000  & 3.81E+00 & 2.52E+07 \\  
                               ~ & 96 & 192000000 & 3.90E+00 & 4.92E+07 \\   \hline
    \end{tabular}
\end{table}

\begin{table}[!ht]
    \caption{Scalability data when all the computations were done on GPUs with MPI-3 
    messaging between CPUs. The speedups on the GPUs relative to the CPU-based information 
    in Tables XV and XVI are also shown.}
    \vspace{2mm}
    \centering
    \begin{tabular}{ c c c c c c }
    \hline
                  &       &          &           & GPU v/s CPU    &  GPU v/s CPU \\
                  &       &          &    GPU    & Speed up       &  Speed up    \\
        A100 GPUs & Zones & GPU time & Zones/sec & Only MPI       &  OpenMP+MPI  \\ \hline
        12 & 24000000  & 5.87E-01 & 4.09E+07 & 4.17 & 5.53 \\ 
        24 & 48000000  & 6.24E-01 & 7.70E+07 & 6.07 & 5.80 \\ 
        48 & 96000000  & 6.70E-01 & 1.43E+08 & 9.63 & 5.69 \\ 
        96 & 192000000 & 7.75E-01 & 2.48E+08 & 14.7 & 5.03 \\ \hline
    \end{tabular}
\end{table}

    Table XV shows the scalability study when only MPI parallelization was used on the
CPUs. Table XV1 shows the same scalability when a hybrid OpenMP+MPI parallelization
paradigm was used on the CPUs; we see that the results in Table XVI are a nice improvement
over Table XV. Table XVII shows the scalability study when only MPI was used across the
CPUs but all the computation was done on GPUs. (As has been explained before, the structure of
the application was such that it was not possible to use MPI natively on the GPUs.) The speed
advantages of the GPUs relative to the two CPU approaches (with different parallelization
strategies) are shown in the last two columns of Table XVII. From the second to last column of
Table XVII we see that using a parallelization model on CPUs that is based exclusively on MPI
is not a winning strategy because the CPU performance degrades with increasing numbers of
CPUs. This is because the patches on the CPUs are small and a lot of ghost zone information has
to be exchanged via a lot of messages. However, from the last column in Table XVII we see that
when a hybrid OpenMP+MPI parallelization strategy is followed for CPU computing then the
patches can be made larger on the CPUs and the number of messages exchanged can be
substantially reduced. From the last column of Table XVII we see that if patches of the same size
are used for GPU computing and for CPU computing, then the GPUs are about five to six times
faster than a comparable number of high-end multicore CPUs. Moreover, Table XVII shows this
speed advantage is sustained over a substantial range in the number of GPUs, or CPUs, that are
used for the scalability study.

\begin{figure}[htp]
    \includegraphics[scale=0.53]{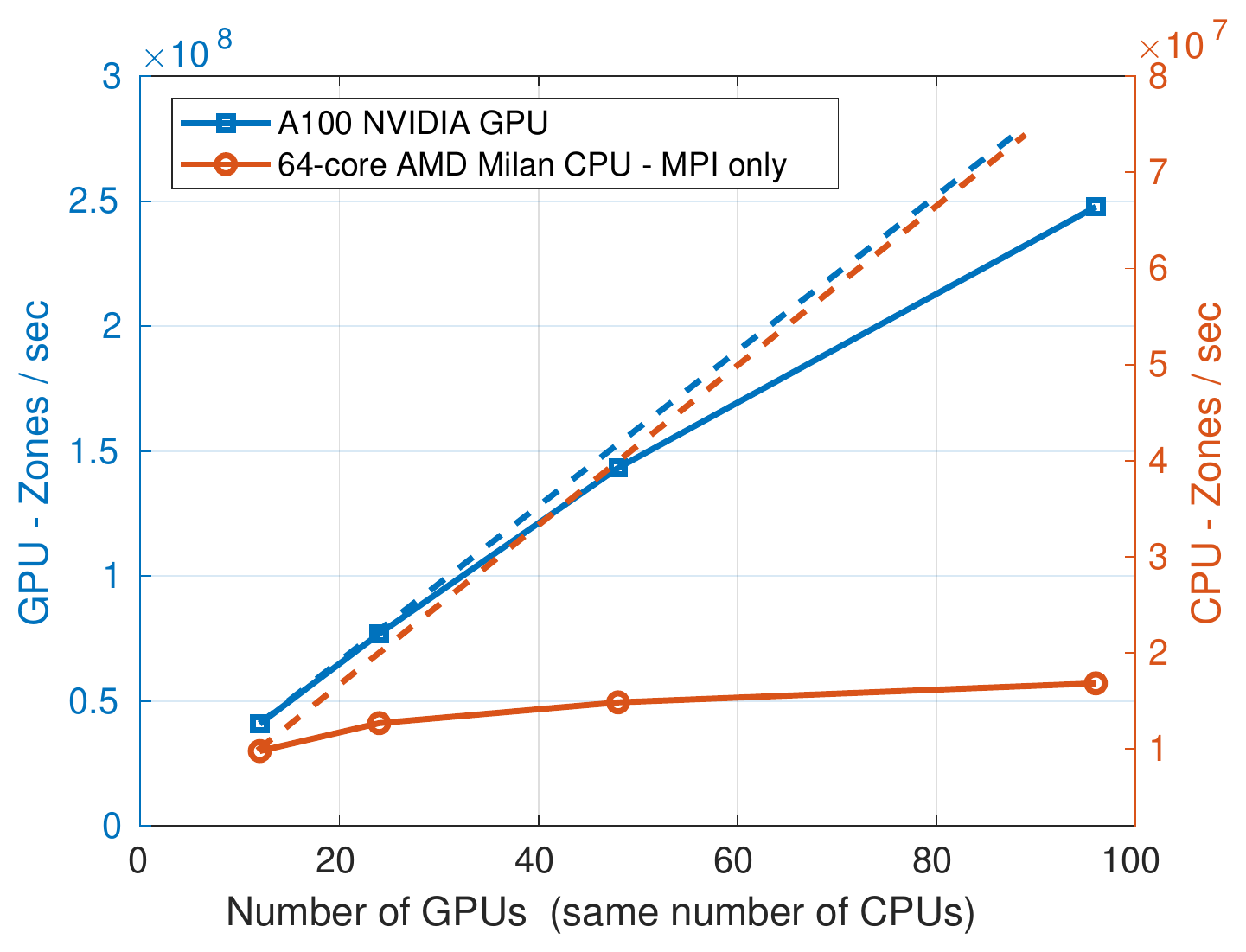} \hspace{4mm}
    \includegraphics[scale=0.53]{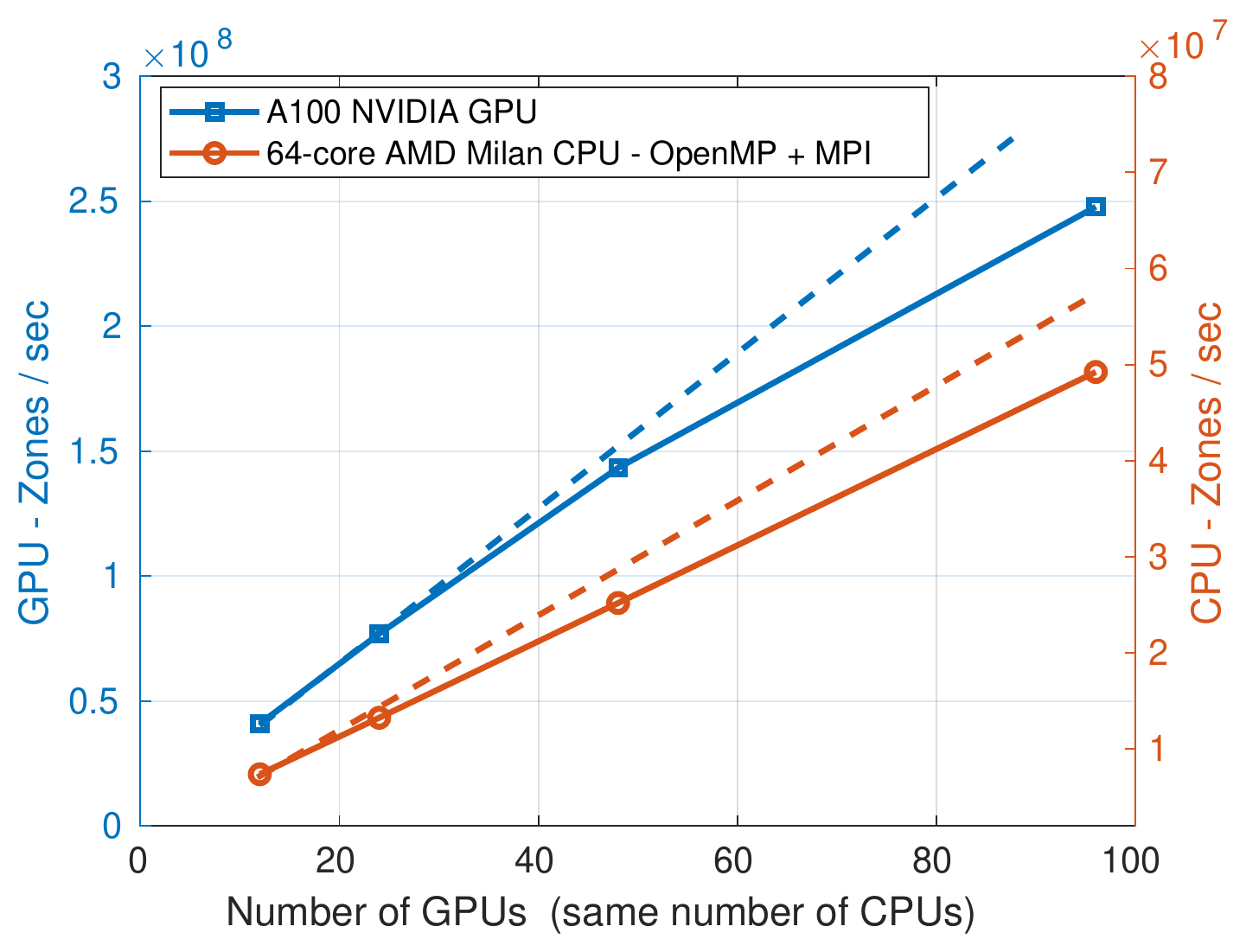}\par
    \hspace{38mm}(a)\hspace{78mm}(b)\par
\caption{Scalability study for the numerical code illustrated in Fig. 8. The scalability study spans 12–96
GPUs. To make this a fair comparison with CPUs, we conducted identical scalability studies that span
12–96 CPUs. The solid blue curves in Figs. 9a and b show the actual GPU results with MPI-3 messaging. 
The dashed blue line shows ideal scalability. In Fig. 9a, the solid red curve shows the actual CPU
results with MPI-3 messaging. Because only MPI was used, the patch sizes were much smaller. The dashed
red line shows the ideal scalability when using MPI-3 messaging. In Fig. 9b, the solid red curve shows
the actual CPU results when a hybrid OpenMP + MPI-3 parallelization model was used. Because we used
OpenMP within the CPU and MPI across CPUs, the patch sizes were much larger. The dashed red line
in Fig. 9b shows the ideal scalability when the hybrid OpenMP + MPI-3 parallelization model was used.
Please note that GPU data should be read off the left vertical axis while CPU data should be read off the
right vertical axis.}
\label{fig:scaling}
\end{figure}

    Figs. 9a and 9b show plots of the information shown in Tables XV, XVI and XVII. The
blue curves in Figs. 9a and 9b show the zones per second updated with the GPU-based runs as a
function of the number of GPUs used. Because the GPU runs are so much faster than the CPU
runs, please use the left vertical axis when viewing GPU data. The red curve in Fig. 9a shows the
zones per second updated as a function of the number of CPUs used when an MPI only model
was used for the parallelism on the CPUs. The red curve in Fig. 9b shows the zones per second
updated as a function of CPUs used when a hybrid OpenMP-MPI model was used for the
parallelism on the CPUs. Because the CPU runs were so much slower than the GPU runs, please
use the right vertical axis when viewing CPU data in Fig. 9. Fig. 9 brings the advantage of GPU
computing, as described in this paper, into very sharp focus!

\section{Conclusions}
    We started with the important realization that GPU computing will be an integral part of
any Exascale supercomputer in the foreseeable future. Since a variety of higher order Godunov
(HOG) schemes are expected to constitute about 30
supercomputers, it is important to explore how such HOG schemes can be ported to such
advanced supercomputer architectures. To that end, we identify three representative areas 
CFD, MHD and CED where HOG schemes are used in quite different styles. We realize that if a
unified pathway can be shown where all such applications can be ported (without too much
effort) to such machines then that would enhance both the value of HOG schemes and their
utilization of GPU-equipped supercomputers.

    To help in the overarching goal of making HOG schemes ready for GPU architectures,
we found that OpenACC is a good vehicle for expressing the parallelism. OpenACC is a
language extension, with the result that a rich set of pragmas enable us to give advice to the GPU
on how to parallelize code. We also acknowledge that successful algorithms should be able to
operate efficiently within the confines of the very small caches that are available on GPUs. The
interconnect between CPU (host) and GPU (device) is also rather slow with the result that tricks
have to be devised to ensure that there is minimal data motion.

    After factoring in the computational landscape within which a HOG scheme has to
operate on Exascale supercomputers, we came to the realization that not every algorithm that is
known to the applied math community will thrive in the constricted environment of GPU
computing. Therefore, it was necessary to identify successful algorithmic choices. WENO
schemes were found to be very favorable approaches for starting with a minimal amount of data
on the GPU and carrying out all of the spatial reconstruction at higher order. Furthermore,
WENO reconstruction can be written in a way so as to optimally work with the very small
caches on GPUs. To avoid frequent data exchanges between host and node and vice versa,
ADER predictor steps were found to be a very favorable way of starting with the spatial
reconstruction and obtaining from it a space-time higher order representation of the PDE. The
Riemann-solver based corrector step had also to be examined and only certain Riemann solvers
(like HLL), which are parsimonious in their use of caches, were found to be favorable.

    Having identified the optimal choices, it is important to show the community how they
can be implemented using the common programming languages that are currently employed.
Therefore, Section II of this paper shows how all parts of a HOG scheme for CFD are easily
implemented using OpenACC extensions to Fortran code. The Appendix does the same for
C/C++ code. With that, this paper serves as a one-stop-shop for giving the reader advice on all
the essential OpenACC pragmas that s/he will need to make a high-quality OpenACC
implementation.

    Just a naïve parallelization with OpenACC will not unlock the dramatic performance
gains that GPUs potentially offer. The most important trick that we found was to identify
skinny variables that move the absolute minimum data from host to node and vice versa.
Without this trick, Table III shows that the GPU speedup invariably suffers. Furthermore, it
degrades with increasing order of accuracy. When this trick is properly used, Table IV shows
that when we go from second to third order we actually see a performance gain! This is because
we are able to utilize GPU cores much more efficiently to carry out the floating point intensive
operations of higher order schemes. GPUs do indeed have very tiny caches (compared to CPU
caches). As a result, Table V identifies a Predictor bottleneck. In other words, while present
generation WENO schemes and Riemann solvers can function efficiently on very small caches,
the ADER predictor algorithm still cannot work very efficiently with the very small caches on a
GPU. We do prove, however, that ADER can be written so that it functions very efficiently with
the larger caches on a CPU. In that sense, the paper also serves to identify the improvements that
future GPU architectures and future algorithmic designs can make to close this gap.

    We have also realized that not everyone might want to go through the intricacies of
implementing an ADER predictor step. For that reason, we also explored the performance of
Runge-Kutta timestepping for CFD and MHD on GPUs. Since this form of timestepping
bypasses the ADER step, we have been able to show that the speedup for third order Runge-
Kutta based schemes is better than the speedup for second order Runge-Kutta based schemes.
Likewise, the speedup for fourth order Runge-Kutta based schemes is better than the speedup for
third order Runge-Kutta based schemes. Overall, however, the ADER timestepping is faster than
the Runge-Kutta timestepping on both CPUs and GPUs. Section V ties all the previous sections
together by describing the process of carrying out a scalability study on a GPU-rich
supercomputer and comparing it to an analogous CPU-based scalability study. This section is
very useful because it shows us how substantial the gains can be if GPUs are mastered, but it
also shows us how one has to dextrously circumvent some of the limitations of GPUs.

    The paper began by saying that GPU computing offered gigantic potential gains, if such
gains could be realized. The results obtained show that with the right algorithmic choices, the
right OpenACC directives and the right tricks, these gains can be mostly realized for a broad
range of HOG schemes. This is especially true for larger meshes where the A100 GPU shows
better than 100X speedup relative to a single core of a CPU. This speedup is certainly much
better than the typical speedups that are often quoted. Indeed, these fantastic gains are only
realized if we pay careful attention to all aspects of the algorithm as well as all the limitations of
present-day GPU architectures. It is, however, very satisfying that such gains can indeed be
realized.

    GPU computing is an evolving enterprise. As other vendors see the commercial benefits
of GPUs, they will design better GPUs that are more dedicated to HPC and knowledge
discovery. All the results here are based on the A100 GPU. The A100 GPU evolved from the
V100 GPU and had several improvements in bandwidth and processing power compared to its
predecessor. It is inevitable that the A100 GPU will give way to even better GPU products.
While we do not have access to the H100 GPU from NVIDIA, early reports seem to indicate that
it has vastly improved bandwidth to main memory. Such improvements will make it possible to
find efficient solution strategies for some of the other algorithms mentioned here, such as the DG
methods and their variants. GPU computing is, after all, an evolving enterprise. The attention to
algorithms that is given in the early Sections of this paper is intended to be evolutionary and
advisory in spirit. As new GPU architectures come online, we will find more opportunities to
bring other higher order algorithms into the fold of GPU computing. The deep analysis of
OpenACC and its usage in these algorithms will, nevertheless, be of value to all such algorithms.

\appendix
\renewcommand\thefigure{\thesection.\arabic{figure}}
\section{OpenACC directives for C/C++ Codes}
\setcounter{figure}{0}
      The OpenACC implementation of a HOG scheme in the C/C++ language is described
here. The overall structure is identical to that of the Fortran code structure described in section
II.2. The reader who is familiar with C/C++ should be able to pick up all the OpenACC insights
from the figures shown in this Appendix. In this Appendix, we also highlight the differences
between Fortran and C/C++ syntax and use of the OpenACC directives.

\begin{figure}
\begin{lstlisting}[language=C, moredelim={[is][commentstyle]{@}{@}},
                               morekeywords={inline}]
#include<stdio.h>
#include<tgmath.h>

#define nx 50
#define ny 50
#define nz 50

// Function declarations

int Initialize( double U[nz+4][ny+4][nx+4][5][5] );

int U_Skinny_to_U( double U_skinny[nz+4][ny+4][nx+4][5],
                   double U[nz+4][ny+4][nx+4][5][5] );

int Boundary_Conditions( double U[nz+4][ny+4][nx+4][5][5] );
int Limiter            ( double U[nz+4][ny+4][nx+4][5][5] );
int Predictor          ( double U[nz+4][ny+4][nx+4][5][5] );

int Make_Flux_X( double      U[nz+4][ny+4][nx+4][5][5],
                 double Flux_X[nz+4][ny+4][nx+5][5]     );

int Make_Flux_Y( double      U[nz+4][ny+4][nx+4][5][5],
                 double Flux_Y[nz+4][ny+5][nx+4][5]     );

int Make_Flux_Z( double      U[nz+4][ny+4][nx+4][5][5],
                 double Flux_Z[nz+5][ny+4][nx+4][5]     );

int Make_dU_dt( double Flux_X[nz+4][ny+4][nx+5][5],
                double Flux_Y[nz+4][ny+5][nx+4][5],
                double Flux_Z[nz+5][ny+4][nx+4][5],
                double dU_dt[nz+4][ny+4][nx+4][5],
                double dt, double dx, double dy, double dz );

int Update_U_Timestep( double U[nz+4][ny+4][nx+4][5][5],
                       double U_skinny[nz+4][ny+4][nx+4][5],
                       double dU_dt[nz+4][ny+4][nx+4][5]
                double *dt_next, double cfl, double dx, double dy, double dz );
\end{lstlisting}
    \caption{showing the header file {\tt hydro.h}, containing the definitions of the
    domain sizes {\tt nx, ny, nz} in the x-, y-, and z-directions respectively; an
    also the function declarations, which are essential for a hydrodynamic simulation.
    This header file is common to all the C/C++ functions and the main function.}
\end{figure}

      Figure A.1 shows the header file containing the preprocessor definitions of the domain
sizes {\tt nx, ny, nz} and the function declarations. The data type and array sizes are specified in the
function interface as required by the C/C++ syntax, which is different from Fortran. It can be
noted that the variable {\tt dt\_next} in the function {\tt Update\_U\_Timestep( )} is passed by reference as a pointer
since this is a scalar output variable. Such special care is not needed in a Fortran subroutine,
since all the variables are passed by reference by default. Similar to the Fortran code, the
function {\tt U\_skinny\_to\_U( )} moves the minimal data from CPU to GPU at the beginning of
every time step.

\begin{figure}
\begin{lstlisting}[language=C, moredelim={[is][commentstyle]{@}{@}},
                               morekeywords={inline}]
#include "hydro.h"
int main()
{
    double U[nz+4][ny+4][nx+4][5][5], U_skinny[nz+4][ny+4][nx+4][5],
       dU_dt[nz+4][ny+4][nx+4][5], Flux_X[nz+4][ny+4][nx+5][5],
      Flux_Y[nz+4][ny+5][nx+4][5], Flux_Z[nz+5][ny+4][nx+4][5];
    double dt, dt_next, cfl, dx, dy, dz;
    const int n_timesteps = 500;
    @#pragma acc data create( U_skinny, U, dU_dt, \
        Flux_X, Flux_Y, Flux_Z, dt, dt_next, cfl, dx, dy, dz )@

    Initialize( U_skinny, dt, cfl, dx, dy, dz ); 
    @#pragma update device( U_skinny, dt, cfl, dx, dy, dz )@

    for( int it = 0; it < n_timesteps; it++ ){
        @#pragma acc update device( U_skinny )@
        U_Skinny_to_U( U_skinny, U );
        Boundary_Conditions( U );
        Limiter( U );
        Predictor( U );
        Make_Flux_X( U, Flux_X );
        Make_Flux_Y( U, Flux_Y );
        Make_Flux_Z( U, Flux_Z );
        Make_dU_dt ( Flux_X, Flux_Y, Flux_Z, dU_dt, dt, dx, dy, dz );
        Update_U_Timestep( U, U_skinny, dU_dt,  &dt_next, cfl, dx, dy, dz );
        @#pragma acc update host( U_skinny )@
        dt = dt_next;
    }
    return 0;
}
\end{lstlisting}
    \caption{showing the structure of a C/C++ version of the hydro code (in pseudocode
    format) at second order with OpenACC extensions to allow it to perform well on
    a GPU. The overall structure is identical to that of the Fortran version. The
    exceptions are, the arrays dimensions starts from {\tt 0} in C and the zone indexing
    is reversed from {\tt ( i, j, k)} of Fortran to {\tt [k][j][i]} of C, to efficiently
    operate on the row-major array storage format in C. The OpenACC pragmas are identical
    to the Fortran version, except the C/C++ version does not need the ending {\tt \#pragma
    acc end} pragma.}
\end{figure}

      Figure A.2. shows the overall structure of a hydro code in C/C++, which is analogous to
the Fortran structure shown in Figure 1. The OpenACC directives have a different structure in
C/C++, they start with {\tt \#pragma acc} and they do not require an {\tt !\$acc end} statement. So, the
analogous form for {\tt !\$acc end data} in Fig. 1 is not present in Fig. A.2.

\begin{figure}
\begin{lstlisting}[language=C, moredelim={[is][commentstyle]{@}{@}},
                               morekeywords={inline}]
#include "hydro.h"
int U_Skinny_to_U( double U_skinny[nz+4][ny+4][nx+4][5], 
                   double        U[nz+4][ny+4][nx+4][5][5] )
{
    int i, j, k, n;

    @#pragma acc parallel vector_length(64) present( U_skinny, U )@
    @#pragma acc loop gang vector collapse(3) independent \
            private( i, j, k, n )@
    for(k = 0; k < nz+4; k++){
        for(j = 0; j < ny+4; j++){
            for(i = 0; i < nx+4; i++){

                @#pragma acc loop seq@
                for(n = 0; n < 5; n++){
                    U[k][j][i][n][0] = U_skinny[k][j][i][n];
                }
            }
        }
    }
    return 0;
}
\end{lstlisting}
    \caption{shows the structure of the {\tt U\_Skinny\_to\_U} function. The array {\tt U\_skinny}
    contains the minimum data that should be carried from the host to the device.
    This subroutine just copies the contents of {\tt U\_skinny} to {\tt U( :, :, :, :, 1)}. The
    OpenACC pragmas are identical to that of the Fortran version.}
\end{figure}

      Figure A.3. shows the C/C++ code for the {\tt U\_skinny\_to\_U( )} function. This is
equivalent to the Fortran code shown in Figure 2. The OpenACC directives that are outside the
triply nested loop are identical in C/C++ and Fortran. The Fortran has the {\tt :} operator for the
whole array operation, whereas the C/C++ does not have any special operator for this. Therefore,
we need another {\tt for} loop inside the triply nested loop. This inner for loop will be executed
serially. It has been marked with the OpenACC directive {\tt \#pragma acc loop seq}.

\begin{figure}
\begin{lstlisting}[language=C, moredelim={[is][commentstyle]{@}{@}},
                               morekeywords={inline}]
#include "hydro.h"
@#pragma acc routine seq@
static inline double MC_limiter( double a, double b, double mc_coef ){
    // Build and return the reconstructed value
}

int Limiter( double U[nz+4][ny+4][nx+4][5][5] )
{
    int i, j, k, n;
    double compression_factor, a, b;

    @#pragma acc data copyin( n, compression_factor )@

    for(n = 0; n < 5; n++){

        if( n == 0 ) compression_factor = 2.0;
        else compression_factor = 1.5;

        @#pragma acc update device( n, compression_factor )@

        @#pragma acc parallel vector_length(64) present( n, U, \
                compression_factor )@
        @#pragma acc loop gang vector collapse(3) independent \
                private( i, j, k, a, b )@
        for(k = 1; k < nz+3; k++){
            for(j = 1; j < ny+3; j++){
                for(i = 1; i < nx+3; i++){
                    a = + U[k][j][i+1][n][0] - U[k][j][i][n][0];
                    b = - U[k][j][i-1][n][0] + U[k][j][i][n][0];
                    U[k][j][i][n][2] = MC_limiter( a, b, compression_factor );
                }
            }
        }
    }
    // Do similarly for other modes in the y- and z-directions
    return 0;
}
\end{lstlisting}
    \caption{shows the C/C++ version of the application of a simple Monotone-Centered
    limiter in the x-direction to all the fluid variables. The inline function {\tt MC\_Limiter}
    is made to accommodate a compression factor which can indeed be changed based
    on the flow field being limited. The OpenACC directives are identical to that
    of the Fortran code, except the C version does not need the ending {\tt \#pragma acc end} 
    pragma. Also, the inline function must be identified as {\tt \#pragma acc routine seq}
    pragma, since the C does not differentiate between a function and a subroutine like Fortran does.}
\end{figure}

      Figure A.4. shows the overall structure of a {\tt Limiter( )} function in C/C++. This is
analogous to the Fortran subroutine shown in Figure 3. The OpenACC directives are identical,
except for the {\tt !\$acc end} directive, which is not needed in C/C++. Unlike Fortran, the C/C++
does not allow functions to be defined inside another function. So, the {\tt MC\_limiter( )} function
is defined as an inline function, outside the {\tt Limiter( )}. The {\tt MC\_limiter( )} is called inside a
triply nested loop serially. So, it is marked with {\tt \#pragma acc routine seq}.

\begin{figure}
\begin{lstlisting}[language=C, moredelim={[is][commentstyle]{@}{@}},
                               morekeywords={inline}]
#include "hydro.h"
@#pragma acc routine seq@
int Predictor_ptwise( double U_ptwise[5][5] ){
    // Predictor step for each zone
}

int Predictor( double U[nz+4][ny+4][nx+4][5][5] )
{
    int i, j, k, n, m;
    double U_ptwise[5][5];

    @#pragma acc routine( Predictor_ptwise ) seq@

    @#pragma acc parallel vector_length(64) present( U )@
    @#pragma acc loop gang vector collapse(3) independent \
            private( i, j, k, n, m, U_ptwise )@
    for(k = 1; k < nz+3; k++){
        for(j = 1; j < ny+3; j++){
            for(i = 1; i < nx+3; i++){
                @#pragma acc loop seq@
                for(n = 0; n < 5; n++)
                    for(m = 0; m < 4; m++)
                        U_ptwise[n][m] = U[k][j][i][n][m];

                Predictor_ptwise( U_ptwise );

                @#pragma acc loop seq@
                for(n = 0; n < 5; n++)
                    U[k][j][i][n][4] = U_ptwise[n][4];
            }
        }
    }
    return 0;
}
\end{lstlisting}
    \caption{shows the C/C++ version of the predictor step for building the 5th mode,
    i.e., the mode that contains the time rate of change. The OpenACC directives
    are identical to that of the Fortran code, except in the C version the inner
    for-loops need to be marked with the {\tt acc loop seq} pragma indicating their sequential
    execution. This is not needed in the Fortran version, where the {\tt :} operator
    is used  for such operations.}
\end{figure}

      Figure A.5. shows the overall structure of a {\tt Predictor( )} function in C/C++. This is
equivalent to the Fortran subroutine, shown in Figure 4. The overall structure and the OpenACC
directives are identical to the Fortran code. As mentioned before, the C/C++ lacks the {\tt :}
operator for array operations. Due to this, {\tt for} loops are used for the array operations, inside the
triply nested loop. These are marked with {\tt \#pragma acc loop seq} to indicate their serial
execution.

\begin{figure}
\begin{lstlisting}[language=C, moredelim={[is][commentstyle]{@}{@}},
                               morekeywords={inline}]
#include "hydro.h"
@#pragma acc routine seq@
int Riemann_Solver_Ptwise( double U_L[5], double U_R[5], 
                                          double flux_x_ptwise[5] ){
    // Build the HLL/HLLEM flux
}

int Make_Flux_X( double U[nz+4][ny+4][nx+4][5][5],
                 double Flux_X[nz+4][ny+4][nx+5][5] )
{
    int i, j, k, n;
    double U_L[5], U_R[5], flux_x_ptwise[5];

    @#pragma acc routine( Riemann_Solver_Ptwise ) seq@

    @#pragma acc parallel vector_length(64) present( U, Flux_X )@
    @#pragma acc loop gang vector collapse(3) independent \
            private( i, j, k, n, U_L, U_R, flux_x_ptwise )@
    for(k = 2; k < nz+2; k++){
        for(j = 2; j < ny+2; j++){
            for(i = 2; i < nx+3; i++){

                @#pragma acc loop seq@
                for(n = 0; n < 5; n++){

                    U_L[n] = U[k][j][i-1][n][0] + 0.5*U[k][j][i-1][n][1] \
                                                + 0.5*U[k][j][i-1][n][4];

                    U_R[n] = U[k][j][i][n][0] - 0.5*U[k][j][i][n][1] \
                                              + 0.5*U[k][j][i][n][4];
                }

                Riemann_Solver_Ptwise( U_L, U_R, flux_x_ptwise );

                @#pragma acc loop seq@
                for(n = 0; n < 5; n++){
                    Flux_X[k][j][i][n] = flux_x_ptwise[n];
                }
            }
        }
    }
    return 0;
}
\end{lstlisting}
    \caption{shows the C/C++ version of the pseudocode for the subroutine to build
    the flux along the x-faces. The {\tt U\_L} variable stores the value of {\tt U} obtained
    from the zone which is at the left to the face and the {\tt U\_R} variable stores
    the value of {\tt U} obtained from the zone which is at the right to the face.  Thus,
    the facial values of the conserved variable {\tt U} are stored in {\tt U\_L, U\_R}. This
    is then fed into the Riemann solver function to build the flux for the corresponding
    face. The OpenACC directives are identical to that of the Fortran code, except
    in the C version the inner for-loops need to be marked with the {\tt acc loop seq}
    pragma indicating their sequential execution. This is not needed in the Fortran
    version, where the {\tt :} operator is used  for such operations. In addition, it
    can be noted that the index of small arrays such as {\tt U\_L, U\_R, flux\_x\_ptwise}
    starts with {\tt 0} instead of {\tt 1}.}
\end{figure}

      Figure A.6. shows the C/C++ code for the {\tt Make\_Flux( )} function. This is equivalent to
the Fortran code shown in Figure 5. The OpenACC directives are identical to the Fortran code.
An additional {\tt \#pragma acc loop seq} is used for the inner loop, which performs the operations
on each fluid component.

\begin{figure}
\begin{lstlisting}[language=C, moredelim={[is][commentstyle]{@}{@}},
                               morekeywords={inline}]
#include "hydro.h"
int Make_dU_dt( double Flux_X[nz+4][ny+4][nx+5][5],
                double Flux_Y[nz+4][ny+5][nx+4][5],
                double Flux_Z[nz+5][ny+4][nx+4][5], 
                double  dU_dt[nz+4][ny+4][nx+4][5], 
                double dt, double dx, double dy, double dz )
{
    int i, j, k, n;

    @#pragma acc parallel vector_length(64) present(      \
            dt, dx, dy,dz, Flux_X, Flux_Y, Flux_Z, dU_dt )@
    @#pragma acc loop gang vector collapse(3) independent \
            private( i, j, k, n )@
    for(k = 2; k < nz+2; k++){
        for(j = 2; j < ny+2; j++){
            for(i = 2; i < nx+2; i++){

                @#pragma acc loop seq@
                for(n = 0; n < 5; n++){
                    dU_dt[k][j][i][n] = \
                     - dt * (Flux_X[k][j][i+1][n] - Flux_X[k][j][i][n]) / dx \
                     - dt * (Flux_Y[k][j+1][i][n] - Flux_Y[k][j][i][n]) / dy \
                     - dt * (Flux_Z[k+1][j][i][n] - Flux_Z[k][j][i][n]) / dz;
                }
            }
        }
    }
    return 0;
}
\end{lstlisting}
    \caption{shows the C/C++ version of the function which builds the {\tt dU\_dt} update
    from the fluxes obtained from the {\tt Make\_Flux\_X, Make\_Flux\_Y and Make\_Flux\_Z}
    functions. The {\tt :} operator of the Fortran is replaced with a sequential inner
    for-loop.}
\end{figure}

      Figure A.7. shows the equivalent C/C++ code for the {\tt Make\_dU\_dt( )} function. This is
analogous to the Fortran shown in Figure 6. The overall structure and the OpenACC directives
are identical between the C/C++ and Fortran.

\begin{figure}
\begin{lstlisting}[language=C, moredelim={[is][commentstyle]{@}{@}},
                               morekeywords={inline}]
#include "hydro.h"
int Update_U_Timestep( double U[nz+4][ny+4][nx+4][5][5],
                double U_skinny[nz+4][ny+4][nx+4][5], 
                double dU_dt[nz+4][ny+4][nx+4][5],
                double *dt_next, double cfl, double dx, double dy, double dz)
{
    int i, j, k, n;
    double U_ptwise[5], dt1 = 1.0e32;
    
    @#pragma acc routine( Eval_Tstep_ptwise ) seq@
    @#pragma acc data copyin( dt1 )@

    @#pragma acc parallel vector_length(64) present(      \
            U, U_skinny, dU_dt, cfl, dx, dy, dz )@
    @#pragma acc loop gang vector collapse(3) independent \
            private( i, j, k, n, U_ptwise ) \
            reduction( min:dt1 )@
    for(k = 2; k < nz+2; k++){
        for(j = 2; j < ny+2; j++){
            for(i = 2; i < nx+2; i++){

                @#pragma acc loop seq@
                for(n = 0; n < 5; n++){
                    U[k][j][i][n][0] += dU_dt[k][j][i][n];
                    U_skinny[k][j][i][n] = U[k][j][i][n][0];
                    U_ptwise[n] = U[k][j][i][n][0];
                }
                dt1 = fmin( dt1, 
                            Eval_Tstep_ptwise( cfl, dx, dy, dz, U_ptwise ) );
            }
        }
    }
    @#pragma acc update host( dt1 )@
    *dt_next = dt1;
    return 0;
}
\end{lstlisting}
    \caption{shows the C/C++ version of the function which updates the {\tt U} from the
    {\tt dU\_dt}. In addition to this, here the {\tt U\_skinny} is updated from the new zone-centered
    values of {\tt U}. The OpenACC directives are identical to that of the Fortran code,
    except in the C version the inner for-loop need to be marked with the {\tt acc loop seq} 
    pragma indicating its sequential execution. Also in this function, the next-timestep
    {\tt dt\_next} is estimated, using the CFL number, zone size and U. The smallest among
    all the zones is found using the {\tt fmin()} function. In order to correctly parallelize
    this reduction operation, the {\tt reduction( min: dt1 )} pragma is used.}
\end{figure}

      Figure A.8. shows the C/C++ code for the {\tt Update\_U\_Timestep( )} function. Here we
perform the time update for the conserved variable {\tt U} and estimate the next timestep,
{\tt dt\_next}. This is analogous to the Fortran code shown in Figure 7. In the C/C++ version, the
scalar variable {\tt dt\_next} is passed by reference as a pointer. In such a case, the function and the
OpenACC directives are unaware of the data structure of the pointer variable. Therefore, to keep
it simple, an additional temporary variable {\tt dt1} is used for the reduction operation. The variable
{\tt dt1} is copied to GPU at the beginning, using the directive {\tt \#pragma acc data copyin( dt1 )}. At
the end of the triply nested loop, the final data stored in {\tt dt1} is copied back to the CPU, using
the OpenACC directive {\tt \#pragma acc update host( dt1 )}. Then it is copied to the {\tt dt\_next}
variable.

\section*{Acknowledgements}
DSB acknowledges support via NSF grants NSF-19-04774, NSF-AST-2009776, NASA-2020-
1241 and NASA grant 80NSSC22K0628. DSB and HK acknowledge support from a Vajra award.



\end{document}